\documentclass[12pt]{elsarticle}
\journal{arXiv and then for publication in impacted journal.}
\usepackage{amsmath, amssymb, enumerate, amsthm,booktabs}

\usepackage{lineno,hyperref}
\modulolinenumbers[5]

\RequirePackage{enumerate}
\RequirePackage{mathrsfs}
\RequirePackage{amsmath,stackrel}
\RequirePackage{amsfonts}
\RequirePackage{multirow}
\RequirePackage{pxfonts}
\RequirePackage{epsfig}
\RequirePackage{color}
\RequirePackage{float}
\RequirePackage{graphics}
\RequirePackage{pgf,tikz}

\newtheorem{Prop}{Proposition}

\RequirePackage{cite}

\begin{document}
\begin{frontmatter}
\title{Relative variation indexes for multivariate continuous distributions on $[0,\infty)^k$ and extensions}
\author[cck]{C\'elestin C. \textsc{Kokonendji}}
\address[cck]{Laboratoire de math\'ematiques de Besan\c con, Universit\'e Bourgogne Franche-Comt\'e, Besan\c con, France}
\ead{celestin.kokonendji@univ-fcomte.fr}
\author[cck]{{Aboubacar Y. \textsc{Tour\'e}}\corref{mycorrespondingauthor}}
\cortext[mycorrespondingauthor]{Corresponding author at: Laboratoire de math\'ematiques de Besan\c con - UMR 6623 CNRS-UFC, UFR Sciences et Techniques - UFC, 16 route de Gray, 25030 Besan\c con Cedex, France; Tel. +33 381 666 341, Fax +33 381 666 623}
\ead{aboubacar\_yacouba.toure@univ-fcomte.fr}
\author[asa]{Amadou \textsc{Sawadogo}}
\address[asa]{UFR de Math\'ematiques et Informatique, Universit\'e  F\'elix Houphou\"et Boigny, 22 BP 582 Abidjan 22, C\^ote d'Ivoire}
\ead{amadou.sawadogo@gmail.com}

\begin{abstract}
We introduce some new indexes to measure the departure of any multivariate continuous distribution on non-negative orthant from a given reference one such the uncorrelated exponential model, similar to the relative Fisher dispersion indexes of multivariate count models. 
The proposed multivariate variation indexes are scalar quantities, defined as ratios of two quadratic forms of the mean vector and the covariance matrix. They can be used to discriminate between continuous positive distributions. Generalized and multiple marginal variation indexes with and without correlation structure, respectively, and their relative extensions are discussed. The asymptotic behavior and other properties are studied. Illustrative examples and numerical applications are analyzed under several scenarios, leading to appropriate choices of  multivariate models. Some concluding remarks and possible extensions are made.
\end{abstract}
\end{frontmatter}


\section{Introduction}\label{sect-Intro}

The choice of a multivariate model from a dataset is not an easy task (e.g., Kotz et al., 2000; Joe, 2014). In practice, we sometimes need simple and effective indicators of multivariate distribution classes in this jungle. They must be appropriate summaries of the multivariate dataset.

Behind the Gaussian distribution and similar to the Poisson distribution for count models (e.g., Kokonendji, 2014), we probably have the exponential distribution on the positive half real line which is the most common probability distributions for this support. It is a particular case of many ones, for instance the lognormal and Weibull distributions, and it has also a wide range of statistical applications in many fields such the reliability; see, e.g., the monograph of Balakrishnan and Basu (1995) for a review. In the multivariate setting, there is not a unique way to define a multivariate exponential distribution; e.g., Basu (1988) and Cuenin et al. (2016).

Recently, Abid et al. (2019abc) have introduced the variation index (VI) for measuring the departure of any absolutely continuous probability distribution concentrated on the non-negative half real line from the equivaried exponential model. Defined as the ratio of variance to squared mean and can be seen as the square of the well-known coefficient of variation (Pearson, 1896), the so-called J\o rgensen variation index (or simply VI) makes it possible to discriminate between univariate continuous distributions to over- and under-variation with respect to exponential distribution and to make inference; see Tour\'e et al. (2019). Since both univariate concepts of VI and of the well-known Fisher (1934) dispersion index with respect to the equidispersed Poisson model are similar (e.g., Tour\'e et al., 2019), we here suggest first a useful and appropriate definition of multivariate over-, equi- and under-variation following the multivariate dispersion indexes of Kokonendji and Puig (2018). Then, we mainly propose an extension for unifying multivariate dispersion and multivariate variation indexes in the framework of natural exponential families.

The rest of the paper is organized as follows. Section \ref{sec:1} presents notations, generalized and relative variation indexes with their interpretation and properties for practical handling. Section \ref{sec:3} illustrates calculations of these measures on some usual bi- and multi-variate continuous positive orthant distributions such beta, exponential and Weibull.  Section \ref{sec:2} provides asymptotic properties of the corresponding estimators. Section \ref{sect-NumAppl} presents example applications from real life and simulated continuous (non-negative  orthant) datasets under several scenarios, and produces some simulation studies. Section \ref{sect-Concl} concludes with some remarks and a unified variability index which includes all multivariate dispersion and variation indexes. 
To make the paper self-contained and more understandable, three appendices are added: (A) a broader multivariate exponential distribution which is derived from Cuenin et al. (2016), (B) a construction of the generalized VI is deduced from Albert and Zhang (2010), and (C) proofs of the asymptotic results are adapted from Kokonendji and Puig (2018).

\section{Multivariate variation indexes}\label{sec:1}

Let $\boldsymbol{Y} = (Y_1,\ldots,Y_k)^{\top}$ be a non-negative continuous $k$-variate random vector on $[0,\infty)^k$, $k\geq 1$. We consider
the following notations: $\sqrt{\mathrm{var}\boldsymbol{Y}}=(\sqrt{\mathrm{var}Y_1},\ldots,\sqrt{\mathrm{var}Y_k})^{\top}$
is the elementwise square root of the variance vector of $\boldsymbol{Y}$; 
$\mathrm{diag}\sqrt{\mathrm{var}\boldsymbol{Y}}=\mathrm{diag}_k (\sqrt{\mathrm{var}Y_j})$
is the $k\times k$ diagonal matrix with diagonal entries $\sqrt{\mathrm{var}Y_j}$ and $0$ elsewhere; and, 
$\mathrm{cov}\boldsymbol{Y}= (\mathrm{cov}(Y_i,Y_j) )_{i,j\in \{1,\ldots,k\}}$
denotes the covariance matrix of $\boldsymbol{Y}$ which is a
$k\times k$ symmetric matrix with entries $\mathrm{cov}(Y_i,Y_j)$ such that
$\mathrm{cov}(Y_i,Y_i)=\mathrm{var}Y_i$ is the variance of $Y_i$. Then
\begin{equation}\label{CovRhoVar}
\mathrm{cov}\boldsymbol{Y}= (\mathrm{diag}\!\sqrt{\mathrm{var}\boldsymbol{Y}})( \boldsymbol{{\textcolor{black}{\rho}}_{\boldsymbol{Y}}})(\mathrm{diag}\!\sqrt{\mathrm{var}\boldsymbol{Y}}),
\end{equation}
where
$\boldsymbol{{\textcolor{black}{\rho}}_{\boldsymbol{Y}}}=\boldsymbol{\rho}(\boldsymbol{Y})$
is the correlation matrix of $\boldsymbol{Y}$; see, e.g., Johnson and Wichern (2007, Eq. 2-36). Note that there are infinitely many multivariate distributions with exponential margins. We denote a generic $k$-variate exponential distribution by $\mathscr{E}_{k}(\boldsymbol{\mu},\boldsymbol{\rho})$, given specific positive mean vector $\boldsymbol{\mu}^{-1}:=(\mu_1^{-1},\ldots,\mu_k^{-1})^{\top}$ and correlation matrix $\boldsymbol{\rho}=(\rho_{ij})_{i,j \in \{1, \ldots, k\}}$; see, e.g., Appendix A for a broader one. The uncorrelated or independent $k$-variate exponential will be written as $\mathscr{E}_k(\boldsymbol{\mu})$, for 
$\boldsymbol{\rho}=\boldsymbol{I}_k$ the $k\times k$ unit matrix.

\subsection{Basic definitions\label{ssec:definition}}

Proceeding along similar lines as Albert and Zhang (2010) and also as Kokonendji and Puig (2018), we define the \textit{generalized variation index} of $\boldsymbol{Y}$ by
\begin{equation}\label{GVI-def}
\mathrm{GVI}(\boldsymbol{Y}) =\frac{\mathbb{E}\boldsymbol{Y}^\top\, ( \mathrm{cov}\boldsymbol{Y})\;\mathbb{E}\boldsymbol{Y}}{(\mathbb{E}\boldsymbol{Y}^{\top}\mathbb{E}\boldsymbol{Y})^2};
\end{equation}
see Appendix B for its construction. 
Remark that when $k=1$, $\mathrm{GVI}$ is the univariate variation index VI (Abid et al., 2019b). The \textit{relative (generalized) variation index} is defined, for two continuous random vectors $\boldsymbol{X}$ and $\boldsymbol{Y}$ on the same support $\mathbb{S}=[0,\infty)^k$ with $\mathbb{E}\boldsymbol{X}=\mathbb{E}\boldsymbol{Y}$ and $\mathrm{GVI}(\boldsymbol{X})>0$, by
\begin{equation}\label{MRVI}
\mathrm{RVI}_{\boldsymbol{X}}(\boldsymbol{Y}): =\frac{\mathrm{GVI}(\boldsymbol{Y})}{\mathrm{GVI}(\boldsymbol{X})}\gtreqqless 1;
\end{equation}
i.e., the \textbf{over-} (\textbf{equi-} and \textbf{under-variation}) of $\boldsymbol{Y}$ compared to $\boldsymbol{X}$, and denoted by $\boldsymbol{Y}\succ \boldsymbol{X}$ ($\boldsymbol{Y}\asymp \boldsymbol{X}$ and $\boldsymbol{Y}\prec \boldsymbol{X}$), is realized if $\mathrm{GVI}(\boldsymbol{Y})>\mathrm{GVI}(\boldsymbol{X})$ ($\mathrm{GVI}(\boldsymbol{Y})=\mathrm{GVI}(\boldsymbol{X})$ and $\mathrm{GVI}(\boldsymbol{Y})<\mathrm{GVI}(\boldsymbol{X})$, respectively).
In the framework of the natural exponential family $F_{\boldsymbol{Y}}$ on $[0,\infty)^k$ (e.g., Chapter 54 in Kotz et al. 2000) , generated by the distribution of $\boldsymbol{Y}$ and characterized by its variance function
$\boldsymbol{m}\mapsto\mathbf{V}_{F_{\boldsymbol{Y}}}(\boldsymbol{m})$, the index $\mathrm{GVI}$ can also be rewritten via $\boldsymbol{m}$ and $\mathbf{V}_{F_{\boldsymbol{Y}}}(\boldsymbol{m})$. As for (\ref{GVI-def}), the ``generalized variation 
function'' defined on the mean domain
$\mathbf{M}_{F_{\boldsymbol{Y}}}\subseteq(0,\infty)^k$ to
$(0,\infty)$ by
\begin{equation}\label{GVI_F(m)}
\boldsymbol{m}\mapsto\mathrm{GVI}_{F_{\boldsymbol{Y}}}(\boldsymbol{m}) =\frac{\boldsymbol{m}^\top \{\mathbf{V}_{F_{\boldsymbol{Y}}}(\boldsymbol{m})\} \boldsymbol{m}}{(\boldsymbol{m}^\top\boldsymbol{m})^2} 
\end{equation}
appears to be very useful through this parameterization.

\subsection{Interpretation and properties\label{ssec:Properties}}

Concerning an interpretation of GVI, we first express the denominator of (\ref{GVI-def}) as
$\mathbb{E}\boldsymbol{Y}^{\top}\mathbb{E}\boldsymbol{Y}=\sqrt{\mathbb{E}\boldsymbol{Y}}^{\top}\, ( \mathrm{diag}\mathbb{E}\boldsymbol{Y}) \sqrt{\mathbb{E}\boldsymbol{Y}}$,
using then (\ref{CovRhoVar}) to rewrite $\mathrm{cov}\boldsymbol{Y}$, obtaining
\begin{equation}\label{GVI-rho}
\mathrm{GVI}(\boldsymbol{Y})=\frac{\{(\mathrm{diag}\!\sqrt{\mathrm{var}\boldsymbol{Y}})\mathbb{E}\boldsymbol{Y}\}^{\top} ( \boldsymbol{{\textcolor{black}{\rho}}_{\boldsymbol{Y}}})\,\{(\mathrm{diag}\!\sqrt{\mathrm{var}\boldsymbol{Y}}) \mathbb{E}\boldsymbol{Y}\} }{\left\lbrace [(\mathrm{diag}\sqrt{\mathbb{E}\boldsymbol{Y}})\!\sqrt{\mathbb{E}\boldsymbol{Y}}]^{\top}(\boldsymbol{{\textcolor{black}{I}}}_k)[ (\mathrm{diag}\sqrt{\mathbb{E}\boldsymbol{Y}})\!\sqrt{\mathbb{E}\boldsymbol{Y}}]\right\rbrace  ^2}.
\end{equation}
From (\ref{GVI-rho}), it is clear that
$\mathrm{GVI}(\boldsymbol{Y})$ makes it possible to compare the full variability of $\boldsymbol{Y}$ (in the numerator) with respect to its expected uncorrelated exponential variability (in the
denominator) which depends only on $\mathbb{E}\boldsymbol{Y}$.

Next, the $\mathrm{GVI}$ index can be considered in itself as a notion of
$k$-variate over-, equi- and under-variation. 

\begin{Prop}\label{Thm-GVI-RVI}
For all positive continuous random vector $\boldsymbol{Y}$ on $(0,\infty)^k$,
$k\geq 1$, and $\boldsymbol{X}\sim\mathscr{E}_k(\boldsymbol{\mu})$
then
$\mathrm{RVI}_{\boldsymbol{X}}(\boldsymbol{Y})=\mathrm{GVI}(\boldsymbol{Y})$.
Furthermore, one has
\begin{equation}\label{GVIexp-rho}
\boldsymbol{Y}\sim\mathscr{E}_k(\boldsymbol{\mu},\boldsymbol{\rho})\quad \Rightarrow\quad \mathrm{GVI}(\boldsymbol{Y})=1+\frac{\boldsymbol{\mu}^{-\top}(\boldsymbol{\rho}-\boldsymbol{I}_k)\boldsymbol{\mu}^{-1}}{(\boldsymbol{\mu}^{-\top}\boldsymbol{\mu}^{-1})^2}\;\gtreqqless 1,
\end{equation}
with $\boldsymbol{\mu}^{-\top}:=(\boldsymbol{\mu}^{-1})^\top=(\mu_1^{-1},\ldots,\mu_k^{-1})$.
\end{Prop}
\bigskip
\noindent
\textbf{Proof}. It is trivial from (\ref{GVI-rho}), with $\mathrm{GVI}(\boldsymbol{Y})=1$ if and only if $\boldsymbol{\rho}=\boldsymbol{I}_k$, i.e., $\boldsymbol{Y}\sim\mathscr{E}_k(\boldsymbol{\mu})$. \hfill $\blacksquare$

\bigskip
From Proposition \ref{Thm-GVI-RVI}, the multivariate exponential model $\boldsymbol{Y}\sim\mathscr{E}_k(\boldsymbol{\mu},\boldsymbol{\rho})$ can be over-, equi-
or under-varied (with respect to the uncorrelated exponential) according to its correlation structure. For instance, if $k=2$ then (\ref{GVIexp-rho}) clearly gives the one-to-one relationship, viz.
$$
\mathrm{GVI}(\boldsymbol{Y})=1+\frac{2\rho\mu_1^{-2}\mu_2^{-2}}{(\mu_1^{-2}+\mu_2^{-2})^2}\;\gtreqqless 1
\quad \Leftrightarrow\quad \rho\;\gtreqqless 0.
$$ 

Finally, if we only want to take into account the variation information coming from the margins,  we can modify GVI by replacing $\mathrm{cov}\boldsymbol{Y}$ in (\ref{GVI-def}) with 
$\mathrm{diag}\,\mathrm{var}\boldsymbol{Y}$, that is $\boldsymbol{\rho}=\boldsymbol{I}_k$ in (\ref{CovRhoVar}), obtaining the ``multiple marginal variation index'', viz.
\begin{equation}\label{MVI-def}
\mathrm{MVI}(\boldsymbol{Y}) = \frac{\mathbb{E}\boldsymbol{Y}^{\top}( \mathrm{diag}\,\mathrm{var}\boldsymbol{Y} )\mathbb{E}\boldsymbol{Y}}{(\mathbb{E}\boldsymbol{Y}^{\top}\mathbb{E}\boldsymbol{Y})^2}
=\sum_{j=1}^k\frac{(\mathbb{E}Y_j)^4}{(\mathbb{E}\boldsymbol{Y}^{\top}\mathbb{E}\boldsymbol{Y})^2} \, \mathrm{VI}(Y_j).
\end{equation}

The expression on the right-hand side of (\ref{MVI-def}) provides a representation of MVI as a weighted average of the univariate variation indexes VI of the components. MVI could be used for exploring profile distributions in multiple positive response regression models (Bonat and J\o rgensen, 2016) or in multivariate continuous time series. Similarly to (\ref{GVI_F(m)}), the corresponding ``multiple marginal variation function'' is defined on the mean domain
$\mathbf{M}_{F_{\boldsymbol{Y}}}\subseteq(0,\infty)^k$ to
$(0,\infty)$ by
\begin{equation}\label{MVI_F(m)}
\boldsymbol{m}\mapsto\mathrm{MVI}_{F_{\boldsymbol{Y}}}(\boldsymbol{m}) =\frac{\boldsymbol{m}^\top \{\mathrm{diag}\,\mathbf{V}_{F_{\boldsymbol{Y}}}(\boldsymbol{m})\} \boldsymbol{m}}{(\boldsymbol{m}^\top\boldsymbol{m})^2}.
\end{equation}
In the same way as (\ref{MRVI}), the relative versions of MVI can be introduced.

\section{Illustrations and comments}\label{sec:3}

We will illustrate our variation indexes with two bivariate models and two families of general $k$-variate ones; it will then be seen how the marginal VIs interplay with the correlation structure in the multivariate variation measures discussed previously. Considering $\boldsymbol{Y}=(Y_1,Y_2)^\top$ and using (\ref{MVI-def}), we explicitly write
\begin{equation}\label{DefGVI-bi}
\mathrm{GVI}(Y_1,Y_2)=\mathrm{MVI}(Y_1,Y_2) +\rho \frac{2(\mathbb{E}Y_1)^2(\mathbb{E}Y_2)^2\sqrt{\mathrm{VI}(Y_1)\mathrm{VI}(Y_2)}}{\{(\mathbb{E}Y_1)^2+(\mathbb{E}Y_2)^2\}^2},
\end{equation}
which points out that GVI is not a weighted average of VIs, as MVI is. Note that $\mathrm{GVI}\gtreqqless\mathrm{MVI}$ accordingly to
$\rho\gtreqqless 0$, with $\mathrm{GVI}=\mathrm{MVI}$ for $\rho=0$. Similar remarks hold for the $k$-variate cases, where the correlation matrix $\boldsymbol{\rho}$ is reduced to $\boldsymbol{I}_k$.

\subsection{Bivariate beta distribution of Arnold and Tony Ng (2011)}

The flexible bivariate beta $\boldsymbol{Y}=(Y_1,Y_2)^\top\sim\mathscr{B}_2(\boldsymbol{\alpha})$ of Arnold and Tony Ng (2011) which exhibits both positive and negative correlation between random variables can be defined as follows. Suppose that $U_1,U_2,V_1,V_2$ and $W$ are independent gamma random variables with common unit scale parameter, i.e., $U_j\sim\mathscr{G}_1(\alpha_j,1)$, $V_j\sim\mathscr{G}_1(\alpha_j',1)$, $j=1,2$ and $W\sim\mathscr{G}_1(\alpha_0,1)$ with $\alpha_j>0$, $\alpha_j'>0$ and $\alpha_0>0$. Then, for $j=1,2$ and $\boldsymbol{\alpha}:=(\alpha_0,\alpha_1, \alpha_2, \alpha_1', \alpha_2')$, one has
$$Y_j:=\frac{U_j+V_j}{U_j+V_1+V_2+W}\:\sim\mathscr{B}_1(\alpha_j+\alpha_j',\alpha_0+\alpha_1'+\alpha_2'-\alpha_j').$$

Since the covariance of $Y_1$ and $Y_2$ cannot be expressed in closed form, it has been numerically shown that the correlation $\rho=\rho(Y_1,Y_2)=\rho(\boldsymbol{\alpha})$ belongs into $[-1,1]$. In fact, from the proposed construction, the positive correlations are obtained when $\alpha_1'=\alpha_2'=0$. For negative correlations, one can consider $\alpha_0=0$ with $\alpha_1$ and $\alpha_2$ fixed, and it will be get closer to $-1$ as $\alpha_1'$ and $\alpha_2'$ get larger. See Arnold and Tony Ng (2011) for more details and connected references.

Thus, for given $\rho(\boldsymbol{\alpha})\in [-1,1]$ of $\boldsymbol{Y}=(Y_1,Y_2)^\top\sim\mathscr{B}_2(\boldsymbol{\alpha})$, the direct calculations of GVI($\boldsymbol{Y}$) through (\ref{DefGVI-bi}) and MVI($\boldsymbol{Y}$) via (\ref{MVI-def}) are obtained from the following  first moments and VIs of the univariate beta random variables $Y_j$, $j=1,2$:
$$\mathbb{E}Y_j=\frac{\alpha_j+\alpha_j'}
{\alpha_0+\alpha_1'+\alpha_2'+\alpha_j}>0$$
and 
\begin{eqnarray*}
\mathrm{VI}(Y_j)&=&\left(\frac{1+\alpha_j+\alpha_j'}{1+\alpha_0+\alpha_1'+\alpha_2'+\alpha_j'}-\mathbb{E}Y_j\right)/\mathbb{E}Y_j\\
&=&\frac{\alpha_0+\alpha_1'+\alpha_2'-\alpha_j'}{(1+\alpha_0+\alpha_1'+\alpha_2'+\alpha_j')(\alpha_j+\alpha_j')}>0.
\end{eqnarray*}
Since $Y_1$ and $Y_2$ are over-, equi- and under-varied, then $\boldsymbol{Y}=(Y_1,Y_2)^\top\sim\mathscr{B}_2(\boldsymbol{\alpha})$ can be too in the bivariate sense and according to the values of $\boldsymbol{\alpha}$. A $k$-variate extension of this bivariate beta distribution is also available for tedious calculations of their multivariate variation indexes.

\subsection{Bivariate Weibull distribution of Teimouri and Gupta (2011)}

Consider the bivariate Weibull $\boldsymbol{Y}=(Y_1,Y_2)^\top\sim\mathscr{W}_2(\alpha_1, \alpha_2, \beta_1, \beta_2, \gamma, \delta)$ of Teimouri and Gupta (2011) built by a copula and having the next density 
\begin{eqnarray*}
f_{\boldsymbol{Y}}(y_1,y_2)&=&\frac{\beta_1\beta_2}{\alpha_1\alpha_2}\left(\frac{y_1}{\alpha_1}\right)^{\beta_1-1}\left(\frac{y_2}{\alpha_2}\right)^{\beta_2-1}\exp\left\{-\left(\frac{y_1}{\alpha_1}\right)^{\beta_1}-\left(\frac{y_2}{\alpha_2}\right)^{\beta_2}\right\}\\
& &\times (1+\delta\exp\left[-(\gamma-1)\left\{\left(\frac{y_1}{\alpha_1}\right)^{\beta_1}+\left(\frac{y_2}{\alpha_2}\right)^{\beta_2}\right\}\right]\\
& &\times \left[\exp\left\{-\left(\frac{y_1}{\alpha_1}\right)^{\beta_1}\right\}-\gamma\right]\left[\exp\left\{-\left(\frac{y_2}{\alpha_2}\right)^{\beta_2}\right\}-\gamma\right]),
\end{eqnarray*}
with $\alpha_j>0$ and $\beta_j>0$, $j=1,2$, $\gamma>1$ and $\delta\in[0,1]$. For $\delta=0$ one gets the uncorrelated bivariate Weibull distribution which depends on both scale parameters $\alpha_j$ and shape parameters $\beta_j$, $j=1,2$.

Since we explicitly have the first, second and product moments of $Y_1$ and $Y_2$ (Teimouri and Gupta, 2011), the calculations of GVI($\boldsymbol{Y}$) using (\ref{DefGVI-bi}) and MVI($\boldsymbol{Y}$) through (\ref{MVI-def}) are derived from the following first moments, VIs and correlation of $Y_j$, $j=1,2$:
$$\mathbb{E}Y_j=\alpha_j\Gamma(1+1/\beta_j),$$
$$\mathrm{VI}(Y_j)=\frac{\Gamma(1+2/\beta_j)}{\{\Gamma(1+1/\beta_j)\}^2}-1$$
depending  only on shape parameter $\beta_j$, and 
\begin{eqnarray*}
\rho=\rho(Y_1,Y_2)&=&\delta\Gamma(1+1/\beta_1)\Gamma(1+1/\beta_2) \\
& &\times [\gamma^{-2-1/\beta_1-1/\beta_2}
-(\gamma+1)\gamma^{-1-1/\beta_1}\left\{\gamma^{-1-1/\beta_2}-(\gamma+1)^{-1-1/\beta_2}\right\} \\
& &-(\gamma+1)\gamma^{-1-1/\beta_2}\left\{\gamma^{-1-1/\beta_1}-(\gamma+1)^{-1-1/\beta_1}\right\} \\
& & + (\gamma+1)^2\left\{\gamma^{-1-1/\beta_1}-(\gamma+1)^{-1-1/\beta_1}\right\}
\left\{\gamma^{-1-1/\beta_2}-(\gamma+1)^{-1-1/\beta_2}\right\}]\\
& &\times\left[\left\{\Gamma(1+2/\beta_1)-\Gamma^2(1+1/\beta_1)\right\} \left\{\Gamma(1+2/\beta_2)-\Gamma^2(1+1/\beta_2)\right\}\right]^{-1/2},
\end{eqnarray*}
where $\Gamma(\cdot)$ is the classical gamma function. The above univariate variation indexes $\mathrm{VI}(Y_j)$, $j=1,2$, satisfy the following equivalences:
\begin{equation}\label{WeibBath}
0<\mathrm{VI}(Y_j)\gtreqqless 1 \stackrel{(i)}{\Longleftrightarrow} 0<\beta_j\Gamma(2/\beta_j)\{\Gamma(1/\beta_j)\}^{-2}\gtreqqless 1
\stackrel{(ii)}{\Longleftrightarrow} 0<\beta_j\lesseqqgtr 1,
\end{equation}
for all $\alpha_j>0$ fixed. Indeed, the first equivalence $(i)$ of (\ref{WeibBath}) is derived from the gamma duplication formula or Legendre's doubling formula - revisited (e.g., Abramowitz and Stegun, 1972; Chap. 6), and the second one $(ii)$ stems from the function study (see, e.g., Table~\ref{Tab_VI_weibull}). Note that conditions (\ref{WeibBath}) on the shape parameter of the univariate Weibull distribution are known, in the opposite sense of $\beta_j$ with respect to $1$, for the failure rate (bathtub) curve in reliability; one can refer to Barlow and Proschan (1981) using the standard coefficient of variation. 
Hence, according to its parameters this bivariate Weibull distribution can be over-, equi- and under-varied with respect to the uncorrelated bivariate exponential distribution.
\begin{table}[!htbp]
\caption{\label{Tab_VI_weibull} Some values of equivalences (\ref{WeibBath}) for the univariate  Weibull distribution.}
\begin{center}
\begin{tabular}{lccccccccc}
\toprule
 $\beta$  & 0.1 & 0.3 & 0.5 & 0.8 & 1 & 2 & 4 & 10 & 100 \\
\toprule
$\beta\Gamma(2/\beta)\{\Gamma(1/\beta)\}^{-2}$ & 92378 & 15.12 & 3 & 1.29 & 1 & 0.64 & 0.54 & 0.5072 & 0.5001 \\
$\mathrm{VI}(\beta)$ & 184755 & 29.24 & 5 & 1.59 & 1 & 0.27 & 0.08 & 0.0145 & 0.0002 \\
\bottomrule
\end{tabular}
\end{center}
\end{table}

\subsection{Multivariate exponential distribution of Marshall and Olkin (1967)}

The $k$-variate exponential $\boldsymbol{Y}=(Y_1,\ldots,Y_k)^\top\sim\mathscr{E}_k(\mu_1, \ldots, \mu_k, \mu_0)$ of Marshall and Olkin (1967) is constructed as follows. Let $X_1,\ldots, X_k$ and $Z$ be univariate exponential random variables with parameters $\mu_1>0,\ldots, \mu_k>0$ and $\mu_0\geq 0$ respectively. Then, by setting $Y_j:=X_j+Z$ for $j=1,\ldots,k$, one easily has $\mathbb{E}Y_j=1/(\mu_j+\mu_0)=\sqrt{\mathrm{var}Y_j}$ and $\mathrm{cov}(Y_j,Y_\ell)=\mu_0/\{(\mu_j+\mu_0)(\mu_\ell+\mu_0)(\mu_j+\mu_\ell+\mu_0)\}$ for all $j\neq\ell$. Note that each correlation $\rho(Y_j,Y_\ell)=\mu_0/(\mu_j+\mu_\ell+\mu_0)$ lies in $[0,1]$ if and only if $\mu_0\geq 0$. 

Thus, using (\ref{GVIexp-rho}) appropriately we obtain 
$$
\mathrm{GVI}(\boldsymbol{Y})=1+\frac{\mu_0\sum_{j=1}^k(\mu_j+\mu_0)^{-1}\{\sum_{\ell\neq j}(\mu_j+\mu_\ell+\mu_0)^{-1}(\mu_\ell+\mu_0)^{-1}\}}{\{(\mu_1+\mu_0)^{-2}+\cdots +(\mu_k+\mu_0)^{-2}\}^2} \geq 1\;\;(\Leftrightarrow\mu_0\geq 0),
$$
and through (\ref{MVI-def}) we easily have
$$
\mathrm{MVI}(\boldsymbol{Y})=
\frac{\sum_{j=1}^k(\mu_j+\mu_0)^{-4}}{\sum_{j=1}^k(\mu_j+\mu_0)^{-4}+2\sum_{1\leq j<\ell\leq 1}(\mu_j+\mu_0)^{-2}(\mu_\ell+\mu_0)^{-2}}<1.
$$

Hence, this multivariate exponential model is always under-varied with respect to the MVI and over- or equi-varied with respect to GVI. If $\mu_0=0$ then this $k$-variate exponential distribution is reduced to $\mathscr{E}_k(\boldsymbol{\mu})$ with $
\mathrm{GVI}(\boldsymbol{Y})=1$. However, the assumption of non-negative correlations between components is sometimes insufficient for some analyzes. We can refer to Appendix A for a more extensive exponential model which is derived as a particular case of a full multivariate Tweedie (1984) models with flexible dependence structure (Cuenin et al., 2016).

\subsection{Multiple stable Tweedie (MST) models}\label{ssect:MST}

Consider the huge  $k$-variate MST class of families $F_{\mathbf{p}}=F_{\mathbf{p}}(m_1,\ldots,m_k,\lambda)$ of models which has been  introducted in Boubacar Ma\"{i}nassara and Kokonendji (2014) to extend the so-called normal stable Tweedie (NST) with $\mathbf{p}=(p_1,0,\ldots,0)$ for $p_1\geq 1$. The normal inverse Gaussian (NIG) model is a common particular case of NST with $p_1=3$ (Barndorff-Nielsen, 1997). This MST class contains infinite subclasses of multivariate distributions among others the gamma-MST with $\mathbf{p}=(2, p_2,\ldots,p_k)$ and inverse Gaussian-MST with $\mathbf{p}=(3, p_2,\ldots,p_k)$. We also have some particular models as the multiple gamma with $\mathbf{p}=(2,\ldots,2)$, multiple inverse Gaussian with $\mathbf{p}=(3,\ldots,3)$ and the gamma-Gaussian with $\mathbf{p}=(2,0,\ldots,0)$ of Casalis (1996). In fact, the MST models are composed by a fixed univariate stable Tweedie (1984) variable having a positive mean domain and random variables that, given the fixed one, are real independent stable Tweedie variables, possibly different, with the same dispersion parameter equal to the fixed component. 

Precisely and for short, within the framework of natural exponential families Kokonendji and Moypemna Sembona (2018) have completely characterized the MST models through their variance functions as follows. Let $\mathbf{p}=(p_1,\ldots,p_k)$ with $p_1\geq 1$ and $p_j\in\{0\}\cup [1,\infty)$ for $j=2,\ldots,k$. Then, the variance function of $F_{\mathbf{p}}$ is given by $\mathbf{V}_{F_{\boldsymbol{p}}}(\boldsymbol{m})=\lambda^{1-p_1}m_1^{p_1-2}\boldsymbol{m}\boldsymbol{m}^\top+\mathrm{diag}\,(0,m_1^{1-p_2}m_2^{p_2},\ldots,m_1^{1-p_k}m_k^{p_k})$ for all $\lambda >0$ and $\boldsymbol{m}=(m_1,\ldots,m_k)^\top\in\mathbf{M}_{F_{\boldsymbol{p}}}=(0,\infty)\times M_{F_{p_2}}\times\cdots\times M_{F_{p_k}}$. Therefore, from (\ref{GVI_F(m)}) and (\ref{MVI_F(m)}), one has 
$$
\mathrm{GVI}_{F_{\boldsymbol{p}}}(\boldsymbol{m}) =\lambda^{1-p_1}m_1^{p_1-2}+\frac{\sum_{j=2}^km_1^{1-p_j}m_j^{p_j+2}}{(\boldsymbol{m}^\top\boldsymbol{m})^2} >0
$$
and 
$$
\mathrm{MVI}_{F_{\boldsymbol{p}}}(\boldsymbol{m}) =\frac{\lambda^{1-p_1}m_1^{p_1-2}\sum_{j=1}^km_j^{4}+\sum_{j=2}^km_1^{1-p_j}m_j^{p_j+2}}{\sum_{j=1}^km_j^{4}+\sum_{1\leq j<\ell\leq k}m_j^2m_\ell^2} >0.
$$

According to different classifications of Kokonendji and Moypemna Sembona (2018), several scenarios occur for $k$-variate over-, equi- and under-variation with respect to GVI and MVI. For instance, let $\lambda=1$ and $p_1=2$ for the exponential-MST subclass then one has $
\mathrm{GVI}_{F_{\boldsymbol{p}}}(\boldsymbol{m})>1$ for all $p_j>1$ with $j=2,\ldots,k$. In order to investigate both indexes  GVI and MVI for $k$-variate (semi-)continuous models on $[0,\infty)^k$, we finally exclude cases $p_1=1$ for the Poisson-MST and also for all $p_j=0$ and $p_j=1$, $j=2,\ldots,k$, related to the normal and Poisson components, respectively. Hence, the NST class is removed from this study.

\section{Estimation and asymptotic properties}\label{sec:2}

Let $\boldsymbol{Y}_1, \ldots, \boldsymbol{Y}_n$ be a random sample from $\boldsymbol{Y}$ with support on $(0,\infty)^k$, where for each $i \in \{ 1, \ldots, n\}$, $\boldsymbol{Y}_i = (Y_{i1},\ldots,Y_{ik})^{\top}$. It is common to consider the empirical versions
\begin{equation}\label{empirical-version}
\overline{\boldsymbol{Y}}_n=\frac{1}{n}\sum_{i=1}^n\boldsymbol{Y}_i = (\overline{Y}_1,\ldots,\overline{Y}_k)^{\top}\quad \mathrm{and}\quad \widehat{\mathrm{cov}\boldsymbol{Y}}=\frac{1}{n-1}\sum_{i=1}^n\boldsymbol{Y}_i\boldsymbol{Y}_i^{\top}-\overline{\boldsymbol{Y}}_n\overline{\boldsymbol{Y}}_n^{\top}
\end{equation}
of the mean vector and covariance matrix of $\boldsymbol{Y}$,
respectively. An estimator of
$\mathrm{GVI}(\boldsymbol{Y})$ directly derived from
(\ref{empirical-version}) is given by
\begin{equation}\label{GVI-estim}
\widehat{\mathrm{GVI}}_n(\boldsymbol{Y})=\frac{\overline{\boldsymbol{Y}}_n^\top\:\widehat{\mathrm{cov}\boldsymbol{Y}}\;\overline{\boldsymbol{Y}}_n}{(\overline{\boldsymbol{Y}}_n^{\top}\overline{\boldsymbol{Y}}_n)^2}.
\end{equation}
Since all the univariate positive continuous variables take positive values, we deduce from Cram\'er (1974, pp. 357-358) that
$\widehat{\mathrm{GVI}}_n(\boldsymbol{Y})$ is an asymptotically
unbiased estimator, i.e., $\mathbb{E}\{\widehat{\mathrm{GVI}}_n(\boldsymbol{Y})\}\approx\mathrm{GVI}(\boldsymbol{Y})$. As for the theoretical variance of $\widehat{\mathrm{GVI}}_n$, we would need at least the moments of fourth order of the components of $\boldsymbol{Y}$.

More interestingly, we establish the following central limit and strong consistency results of $\widehat{\mathrm{GVI}}_n$ and $\widehat{\mathrm{MVI}}_n$. The proofs are given in Appendix C.

\begin{Prop}\label{Prop-DeltaGVI}
Let $\boldsymbol{Y} = (Y_1,\ldots,Y_k)^{\top}$ be a positive continuous $k$-variate random vector on $(0,\infty)^k$, $k\geq 1$, such that $\mathbb{E}(Y_{\ell_1}Y_{\ell_2}Y_{\ell_3}Y_{\ell_4})<\infty$. Let also $\boldsymbol{Y}_1,\ldots,\boldsymbol{Y}_n$ be a random sample from $\boldsymbol{Y}$.
\begin{itemize}
\item[(i)]As $n \to \infty$,
$$
\sqrt{n} \, \{\widehat{\mathrm{GVI}}_n(\boldsymbol{Y})-\mathrm{GVI}(\boldsymbol{Y})\} \rightsquigarrow \mathcal{N}(0,\sigma_{gvi}^2),
$$
where $\rightsquigarrow$ stands for convergence in distribution and $\mathcal{N}(0,\sigma_{gvi}^2)$
is the centered normal distribution with variance
$\sigma_{gvi}^2 = \boldsymbol{\Delta}^\top\boldsymbol{\Gamma}\boldsymbol{\Delta}$.
The $\{k+k(k+1)/2\}\times 1$ vector
$\boldsymbol{\Delta}=(\ldots,\Delta_j,\ldots;\ldots,\Delta_{j\ell},\ldots)_{j \in \{1,\ldots,k\};
\ell \in \{j,\ldots,k\}}^\top$ and the $\{k+k(k+1)/2\}\times \{k+k(k+1)/2\}$
four-block symmetric matrix
$$
\boldsymbol{\Gamma}=
\begin{bmatrix}
\boldsymbol{\Sigma} & \boldsymbol{\Gamma}_3^{\top} \\
\boldsymbol{\Gamma}_3 & \boldsymbol{\Gamma}_4
\end{bmatrix}
$$ 
are such that, for all $j,j',j'' \in \{ 1,\ldots,k\}$,$$
\Delta_j = \left\{ 2\sum_{j'=1}^k\mathbb{E}Y_{j'} \, \mathrm{cov}(Y_j,Y_{j'})-4\mathbb{E}Y_{j}\left( \sum_{j'=1}^k(\mathbb{E}Y_{j'})^2\right) \,\mathrm{GVI}(\boldsymbol{Y})\right\}/(\mathbb{E}\boldsymbol{Y}^{\top}\mathbb{E}\boldsymbol{Y} )^2,
$$
$\Delta_{jj}=(\mathbb{E}Y_j)^2/(\mathbb{E}\boldsymbol{Y}^{\top}\mathbb{E}\boldsymbol{Y})^2$,
$\Delta_{j\ell}=2{\mathbb{E}Y_j\mathbb{E}Y_\ell}/(\mathbb{E}\boldsymbol{Y}^{\top}\mathbb{E}\boldsymbol{Y})^2$
for $\ell \in \{j+1,\ldots,k\}$,
$\boldsymbol{\Sigma}(j;j')=\mathrm{cov}(Y_j,Y_{j'})$,
$\boldsymbol{\Gamma}_3^\top(j;j',\ell')=\mathrm{cov}(Y_j,Y_{j'}Y_{\ell'})$
for $\ell' \in \{ j',\ldots,k\}$ and
$\boldsymbol{\Gamma}_4(j',\ell';j'',\ell'')=\mathrm{cov}(Y_{j'}Y_{\ell'},
Y_{j''}Y_{\ell''})$ for $\ell'' \in \{ j'',\ldots,k\}$;
\item [(ii)] As $n \to \infty$,
$$
\sqrt{n} \, \{ \widehat{\mathrm{MVI}}_n(\boldsymbol{Y})-\mathrm{MVI}(\boldsymbol{Y})\} \rightsquigarrow \mathcal{N}(0,\sigma_{mvi}^2),
$$
with
$\sigma_{mvi}^2 = \boldsymbol{\Lambda}^\top\boldsymbol{\Pi}\,\boldsymbol{\Lambda}$.
The $2k\times 1$ vector
$\boldsymbol{\Lambda}=(\Lambda_1,\ldots,\Lambda_k;\Lambda_{11},\ldots,
\Lambda_{kk})^\top$ and the $2k\times 2k$ four-block symmetric matrix
$$\boldsymbol{\Pi}=
\begin{bmatrix}
\boldsymbol{\Sigma} & \boldsymbol{\Pi}_3^{\top} \\
\boldsymbol{\Pi}_3 & \boldsymbol{\Pi}_4
\end{bmatrix}
$$ 
are such that, for $j, j' \in \{1,\ldots,k\}$,
$$\Lambda_j= \left\lbrace  2\mathbb{E}Y_j\mathrm{var}Y_j-4\mathbb{E}Y_j\left( \sum_{j'=1}^k(\mathbb{E}Y_{j'})^2\right)\,\mathrm{MVI}(\boldsymbol{Y})\right\rbrace  /(\mathbb{E}\boldsymbol{Y}^{\top}\mathbb{E}\boldsymbol{Y})^2,$$ 
$\Lambda_{jj}=(\mathbb{E}Y_j)^2/(\mathbb{E}\boldsymbol{Y}^{\top}\mathbb{E}\boldsymbol{Y})^2$,
$\boldsymbol{\Sigma}(j;j')=\mathrm{cov}(Y_j,Y_{j'})$,
$\boldsymbol{\Pi}_3^\top(j;j')=\mathrm{cov}(Y_j,Y_{j'}^2)$ and
$\boldsymbol{\Pi}_4(j;j')=\mathrm{cov}(Y_j^2, Y_{j'}^2)$.
\end{itemize}
\end{Prop}

Note that Parts (i) and (ii) of Proposition \ref{Prop-DeltaGVI} provide the same result for $k=1$ with 
\begin{eqnarray*}
\sigma_{gvi}^2=\sigma_{mvi}^2 & = & \frac{1}{(\mathbb{E}Y)^6} \, [\mathbb{E}Y^4(\mathbb{E}Y)^2-4\mathrm{var}Y\mathbb{E}Y\mathbb{E}Y^3\\ 
& & + \;4\mathrm{var}Y\mathbb{E}Y^2(\mathbb{E}Y)^2-\{\mathbb{E}Y(\mathbb{E}Y^2)\}^2+4(\mathrm{var}Y)^3];
\end{eqnarray*}
see Tour\'e et al. (2019, Part (i) of Section 4.1 with $\lambda=1$. Also, an asymptotic confidence interval for $\mathrm{GVI}(\boldsymbol{Y})$ is
expressed as
$$\left(\widehat{\mathrm{GVI}}_n-u_{1-\alpha/2}\widehat{\sigma}_{gvi}/\sqrt{n},\widehat{\mathrm{GVI}}_n+u_{1-\alpha/2}\widehat{\sigma}_{gvi}/\sqrt{n}\right),
$$
where $u_p$ is the $p$th percentile of the standard normal distribution $\mathcal{N}(0,1)$ and
$\widehat{\sigma}_{gvi}^2=\widehat{\boldsymbol{\Delta}}_n^\top\widehat{\boldsymbol{\Gamma}}_n\widehat{\boldsymbol{\Delta}}_n$
is the corresponding empirical version of $\sigma_{gvi}^2$ (Proposition
\ref{Prop-DeltaGVI}). A similar result also holds for the intuitive index MVI.
Finally, we state the following results for strong consistency.

\begin{Prop}\label{Prop-sConsGVI}
Let $\boldsymbol{Y} = (Y_1,\ldots,Y_k)^{\top}$ be a positive continuous $k$-variate random vector on $(0,\infty)^k$, $k\geq 1$, such that $\mathbb{E}(Y_{\ell_1}Y_{\ell_2})<\infty$. If $\boldsymbol{Y}_1,\ldots,\boldsymbol{Y}_n$ be a random sample from $\boldsymbol{Y}$, then
$$
\widehat{\mathrm{GVI}}_n(\boldsymbol{Y})\stackrel{a.s.}{\longrightarrow} \mathrm{GVI}(\boldsymbol{Y})\;\mathrm{and}\;\widehat{\mathrm{MVI}}_n(\boldsymbol{Y})\stackrel{a.s.}{\longrightarrow} \mathrm{MVI}(\boldsymbol{Y}),\;\mathrm{as}\;n\to\infty,
$$
where $\stackrel{a.s.}{\longrightarrow}$ stands for almost sure convergence.
\end{Prop}

In finite samples, we suggest the use of the bootstrap method for approximating the variance of all the estimators and their corresponding confidence intervals lenght.
\section{Numerical applications}
\label{sect-NumAppl}

All computations have been done with the \textsf{Python} (Python Software Foundation, 2019) and \textsf{R} software (\textsf{R} Core Team, 2018). To generate a $k$-variate continuous positive orthant distribution given $k$ (over-, equi- and under-varied) marginals and correlation matrix $\boldsymbol{\rho}$, we have used the NORmal To Anything (NORTA) method (e.g., Su, 2015). 

In practical way, we will consider the exponential, lognormal and Weibull distribution (e.g., Dey and Kundu, 2009). They have been used quite effectively in analyzing positively skewed data, which play important roles in the reliability analysis. Recall here that the univariate exponential distribution $\mathscr{E}_1(\mu)$ is always equi-varied for all $\mu>0$ and, both univariate lognormal $\mathscr{L\!\!N}_1(m,\sigma^2)$ and Weibull $\mathscr{W}_1(\alpha,\beta)$ models are over- (equi- and under-) varied for $0<\sigma^2\gtreqqless\log 2$ with $m\in\mathbb{R}$ and, from (\ref{WeibBath}) for $0<\beta\lesseqqgtr 1$ with $\alpha >0$, respectively. All these theoretical behaviors work well on simulated datasets of univariate exponential, lognormal and Weibull distributions that we omit presenting here.

\subsection{Some scenarios of bivariate cases and a real $4$-variate dataset}

We first consider Table~\ref{Tab_BiVIsimul} consisting of fifteen simulated bivariate datasets from exponential, lognormal and Weibull distributions, presenting several scenarios of correlation (positive or negative) and marginal over-, equi- or under-variation. The table presents a summary
of these datasets, along with the sample values of the indexes GVI and MVI.

\begin{table}[!htbp]
\caption{\label{Tab_BiVIsimul} Some scenarios from fifteen simulated bivariate datasets with $\widehat{\mathrm{VI}}_j=\widehat{\sigma}_j^2/\bar{y}_j^2$ and marginal variations (MV): Over- (O), Equi- (E) and Under-variation (U).}
\begin{center}
\begin{tabular}{lrrcrcccrccc}
\toprule
  Dataset & $n\;\;$ & $\bar{y}_1$ &$\bar{y}_2$&$\widehat{\sigma}_1^2\:$&$\widehat{\sigma}_{2}^2\:$&$\widehat{\mathrm{VI}}_1\;$&$\widehat{\mathrm{VI}}_2\;$ & $\widehat{\rho}_{12}\;$&MV & $\widehat{\mathrm{MVI}}$ & $\widehat{\mathrm{GVI}}$  \\
\toprule
No~1 &50 &2.58&3.80&13.36&31.08&2.02&2.15 &$-0.23$&O/O & 2.11 &\textbf{1.00} \\
No~2 & 80 &1.34&1.31&3.98&4.05&2.21&2.36 &$-0.33$ &O/O & 2.14&\textbf{0.76} \\
No~3 &100 & 1.94&1.96&3.74&3.76&0.99&0.98 &0.20&E/E & 0.99&\textbf{1.98} \\
No~4 &120 &5.76&5.04&34.97&29.47&1.02&1.01 &$-0.40$&E/E &1.00 &\textbf{0.34} \\
No~5 &150&0.27&0.26&0.02&0.04&0.24&0.51 &0.56&U/U & 0.39 &\textbf{2.28} \\
No~6 &300 &8.77&8.55&22.87&20.82&0.30&0.28 &0.03&U/U & 0.25 &\textbf{0.15} \\
No~7 &500 &5.60&2.10&299.03&20.32&9.53&4.63 &$-0.06$&O/O & 7.40 &\textbf{7.32} \\
No~8 &800 &10.59&10.39&110.44&111.03&0.98&1.02 &0.05&E/E & 0.98 &\textbf{1.01} \\
No~9 &1000 &2.66&1.76&2.35&0.41&0.33&0.13 &0.78&U/U & 0.17 &\textbf{1.02}\\
No~10 & 1000 &6.10 &1.74 &323.00 & 0.52 &8.68&0.17&0.49 & O/U & 7.42 & \textbf{7.50} \\
No~11 &1500 &1.77&3.41&0.44&17.88&0.14&1.55 &$-0.57$&U/O & 0.96 &\textbf{0.87} \\
No~12 &3000 &0.75&0.62&1.62&0.11&2.86&0.28 &$-0.29$&O/U & 1.06 &\textbf{0.99} \\
No~13 &3000 &1.00&1.09&1.00&3.95&1.00&3.33 &$-0.29$&E/O & 1.19 &\textbf{0.82} \\
No~14 &5000 &0.68&1.00&0.42&1.02&0.89&1.02 &0.80&U/E & 0.98 &\textbf{2.23}\\
No~15 & 8000 &1.98 &3.37 &3.83 & 17.97 &0.98&1.58&$-0.38$ & E/O & 1.33 & \textbf{0.99} \\
\bottomrule
\end{tabular}
\end{center}
\end{table}

In order to measure the departure from the bivariate uncorrelated exponential of the considered datasets our estimated index $\widehat{\mathrm{GVI}}$ provides a very good summary of the bivariate variation by taking into account both marginal variation and the non-null correlation value $\widehat{\rho}_{12}$. Indeed, the bivariate equi-variation ($\widehat{\mathrm{GVI}}=1$) is significantly obtained here for both over-varied marginals with negative correlation (No~1), for both under-varied marginals with large positive correlation (No~9), for both equi-varied marginals with very weak   correlation (No~8), and either for one marginal over-varied and the other under-varied (No 12) or equi-varied (No~15) with negative correlations. The bivariate over-variation ($\widehat{\mathrm{GVI}}>1$) is pointed out for both over-varied marginals with weak negative correlation (No~7), for both equi-varied  marginals with positive correlation (No~3), for both under-varied marginals with positive correlation (No~5), and either for one marginal under-varied and the other over-varied (No~10) or equi-varied (No~14) with positive correlation. Concerning the bivariate under-variation ($\widehat{\mathrm{GVI}} \in (0,1)$), this is pointed out for both over-varied marginals with negative correlation (No~2), for both equi-varied  marginals with negative correlation (No~4), for both under-varied marginals with weak positive correlation (No~6), and either for one marginal over-varied and the other under-varied (No~11) or equi-varied (No~13) with negative correlation.
In the common sense, we always have the bivariate over-/under-variation for both over-/under-varied marginals with positive/negative correlation. The values of $\widehat{\mathrm{GVI}}$ provide the corresponding degree of (over-/under-) variation with respect to the reference value~1 of the bivariate equi-variation. For instance, we detect a higher degree of over-variation in No~7 and No~10 than in No~3, No~5 or No~14; similarly, we detect a weaker degree of bivariate under-variation (close to 1) in No~11 and No~13 than in Nos.~4 or 6.

Similarly, the marginal index $\widehat{\mathrm{MVI}}$ also works very well,
summarizing both marginal variations (without correlation). Both indexes $\widehat{\mathrm{MVI}}$ and $\widehat{\mathrm{GVI}}$ are close when the correlation is quasi-null (Nos.~7 or ~8). For the sake of brevity, we omit here an analysis of the standard errors of the estimated indexes of these datasets; a complete analysis will be done in the next section for $6$-variate datasets. 

In summary, multivariate variation indexes MVI and GVI are meaningful because they summarize the variation behavior from each individual variable. In addition, GVI also contains information
about their correlation. They can be used for descriptive analysis, for clustering, for comparing different datasets and for testing departures from known multivariate distributions as Tour\'e et al. (2019) for univariate case.

Secondly, we consider the real $4$-variate dataset 
which refers to the annual observations from 1900 to 1989 of the United Stated. It is reported by Hayashi (2000): the first variable $(m)$ is the  natural log of the money M1, the second $(p)$ is the natural log of the net national product price deflator, the third $(y)$ natural log of the net national product and the fourth $(r)$ is the commercial paper rate in percent at an annual rate.

To measure the departure from the 4-variate uncorrelated exponential distribution of the considered dataset, our estimated indexes 
provide very good summaries through $\widehat{\mathrm{GVI}}=0.1397$ and $\widehat{\mathrm{MVI}}=0.0771$. Indeed, both indexes strongly show a 4-variate under-variation with $0<\widehat{\mathrm{MVI}}<\widehat{\mathrm{GVI}}<1$. Since $\widehat{\mathrm{MVI}}$ is very close to $0$ than $\widehat{\mathrm{GVI}}$, each of the four marginal distributions must be univariate under-varied with the correlation matrice having only  positive coefficients. Table~\ref{SummuryRsd} confirms this analysis only from results of $\mathrm{MVI}$ and $\mathrm{GVI}$. Thus, one can choose an appropriate theoretical 4-variate distribution for modelling this dataset and their (interest) parameters adjust directly by estimation.
\begin{table}[!htbp]
\caption{\label{SummuryRsd} Summary of real $4$-variate data with $\widehat{\mathrm{VI}}_j=\widehat{\sigma}_j^2/\bar{y}_j^2$ and the marginal variation (MV$_j$): Under-variation (U), size $n=90$, $\det\widehat{\boldsymbol{\rho}}=0.0003$,  $\widehat{\mathrm{GVI}}=0.1397$ and
$\widehat{\mathrm{MVI}}=0.0771$.}
\begin{center}
\begin{tabular}{lcccrrrr}
\toprule
  $j$ & $\bar{y}_j$ & $\widehat{\sigma}_j^2$ & $\widehat{\mathrm{VI}}_j$ (MV$_j$) & $\widehat{\rho}_{j1}\;\;$ & $\widehat{\rho}_{j2}\;\;$& $\widehat{\rho}_{j3}\;\;$& $\widehat{\rho}_{j4}\;$\\
\toprule
1 & 4.1476 & 1.9630 & 0.1141 (U) & 1.0000 & 0.9579 & 0.9905 & 0.3926 \\ 
  2 & 3.1709 & 0.6049 & 0.0602 (U) & 0.9579 & 1.0000 & 0.9552 & 0.6002 \\ 
  3 & 2.2610 & 0.6330 & 0.1238 (U) & 0.9905 & 0.9552 & 1.0000 & 0.4331 \\ 
  4 & 4.5547 & 8.4074 & 0.4053 (U) & 0.3926 & 0.6002 & 0.4331 & 1.0000 \\
\bottomrule
\end{tabular}
\end{center}
\end{table}

\subsection{Other multivariate cases and simulation studies}

In this section we first study a $6$-variate simulated dataset. We then analyze the behavior of the asymptotic variances and confidence intervals by simulation. Finally, we compare the asymptotic standard errors of GVI and MVI to those obtained from the bootstrap method.

The $6$-variate dataset of size $n=560$ is simulated following this scenario. We have considered two over-, two equi- and two under-variations as univariate marginals with the  theoretical correlation matrix such that
$$
\det\boldsymbol{\rho} := \det\left(
\begin{array}{cccccc}
 1 & -0.03 & 0.57 & 0.12 & 0.24&-0.57   \\
  -0.03& 1 & -0.05 & -0.09 & 0.36&0.04  \\
 0.57 & -0.05 & 1 & 0.35 & 0.27&-0.58\\
 0.12 & -0.09 & 0.35 & 1 & 0.18&0.04  \\
  0.24 & 0.36 & 0.27 & 0.18 & 1&0.06 \\
  -0.57&0.04&-0.58&0.04&0.06&1\\
 \end{array}
\right)=0.2051.
$$
Table~\ref{SummurySsd}
shows the summary needed to compute the variation indexes $\widehat{\mathrm{GVI}}$ and $\widehat{\mathrm{MVI}}$. As commented before for Table~\ref{Tab_BiVIsimul}, we also observe a different behavior of the two variation indexes in this $6$-variate example. We obtain here
$\widehat{\mathrm{GVI}}=1.0572\approx 1$ and
$\widehat{\mathrm{MVI}}=0.9637\approx 1$, both indicating a $6$-variate phenomenon of quasi-equi-variation. 
\begin{table}[!htbp]
\caption{\label{SummurySsd} Summary of simulated $6$-variate data with marginal variations (MV): Over- (O), Equi- (E) and Under-variation (U), size $n=560$ and $\det\widehat{\boldsymbol{\rho}}=0.2063$ such that $\widehat{\mathrm{GVI}}=1.0572\approx 1$ and
$\widehat{\mathrm{MVI}}=0.9637\approx 1$.}
\begin{center}
\begin{tabular}{lcccrrrrrr}
\toprule
  $j$ & $\bar{y}_j$ & $\widehat{\sigma}_j^2$ & MV$_j$ & $\widehat{\rho}_{j1}\;\;$ & $\widehat{\rho}_{j2}\;\;$& $\widehat{\rho}_{j3}\;\;$& $\widehat{\rho}_{j4}\;$& $\widehat{\rho}_{j5}\;$&$\widehat{\rho}_{j6}\;$ \\
\toprule
1 & 1.2245 & 3.2031 &O&1.0000   &$-0.0197$ & 0.5572    & $0.1074$   & 0.2939 & $-0.5586$\\
2 & 0.4929 & 0.2324 &E&$-0.0197$&1.0000    & $-0.0683$ & $-0.1078$  & 0.3293 &0.0102\\
3 & 1.1548 & 0.4834 &U&0.5572   &$-0.0683$ & 1.0000    & 0.3136     & 0.3116 & $-0.5946$\\
4 & 0.9507 & 0.9236 &E&0.1074   &$-0.1078$ & 0.3136    & 1.0000     & 0.1451 & 0.0264 \\
5 & 4.3871 & 28.6346&O&0.2939   &0.3293    & 0.3116    & 0.1451     & 1.0000 &0.0310 \\
6 & 0.9093 & 0.1039 &U&$-0.5586$& 0.0102   & $-0.5946$ & 0.0264     & 0.0310 & 1.0000 \\
\bottomrule
\end{tabular}
\end{center}
\end{table}
\begin{table}[!htbp]
\caption{\label{SummurySCI} Asymptotic variances and confidence intervals ($u = u_{0.975}=1.96$) from subsamples of simulated $6$-variate data with $n=\mbox{10 000}$ having the same parameters as for Table~\ref{SummurySsd}.}
\begin{center}
\begin{tabular}{rcccccc}
\toprule
  $n\;$ & $\det\widehat{\boldsymbol{\rho}}$ & $\widehat{\sigma}_{gvi}^2$ & $\widehat{\sigma}_{mvi}^2$ & $\widehat{\mathrm{GVI}}_n\pm u\widehat{\sigma}_{gvi}/\!\sqrt{n}$ & $\widehat{\mathrm{MVI}}_n\pm u\widehat{\sigma}_{mvi}/\!\sqrt{n}$ \\
\toprule
50 & 0.2030 & 19163.96 & 19056.40 & 1.2533  $\pm$ 38.3712 & 1.2395 $\pm$ 38.2634 \\
100 & 0.1571 & 36433.89 & 36299.33 & 0.9571 $\pm$ 37.4111 & 0.9189 $\pm$37.3420 \\
300 & 0.2092 & 28413.48 & 28229.35 & 1.2442 $\pm$ 19.0743 & 1.1788 $\pm$19.0124 \\
500 & 0.2050 & 21789.08 & 21618.11 & 1.1487 $\pm$ 12.9385 & 1.0753 $\pm$12.8876 \\
1 000 & 0.1958 & 14589.39 & 14448.38&1.0147 $\pm$ 7.4863 & 0.9366 $\pm$7.4500 \\
3 000 & 0.2017 & 17982.94 & 17800.34&1.1648 $\pm$ 4.7986 & 1.0888 $\pm$4.7742\\
5 000 & 0.2067 & 18892.96 & 18688.16&1.2188 $\pm$ 3.8099 & 1.1316 $\pm$3.7892 \\
10 000 & 0.2067 & 17558.26 & 17354.58&1.1714$\pm$ 2.5971 & 1.0818 $\pm$2.5820 \\
\bottomrule
\end{tabular}
\end{center}
\end{table}
\begin{table}[!htbp]
\caption{\label{SummurySCI2} Asymptotic variances and confidence
intervals ($u = u_{0.975}=1.96$) from subsamples of simulated
over-varied $4$-variate data with $n=\mbox{10 000}$.}
\begin{center}
\begin{tabular}{rccccc}
\toprule
  $n\;$ & $\det\widehat{\boldsymbol{\rho}}$ & $\widehat{\sigma}_{gvi}^2$ & $\widehat{\sigma}_{mvi}^2$ & $\widehat{\mathrm{GVI}}_n\pm u\widehat{\sigma}_{gvi}/\!\sqrt{n}$ & $\widehat{\mathrm{MVI}}_n\pm u\widehat{\sigma}_{mvi}/\!\sqrt{n}$ \\
\toprule
$50$  &0.3209&  6524.25&  3150.98 & 4.1477 $\pm$ 22.3887& 2.4510 $\pm$ 15.5592 \\
$100$ &0.3452&  1697.28&  1190.33& 3.1632 $\pm$ 8.0747& 1.8721 $\pm$ 6.7621 \\
$300$ &0.6915&  5014.56&  3877.69& 3.5238 $\pm$ 8.0132& 2.7728 $\pm$ 7.0465\\
$500$ &0.7071&  6547.90&  1803.12& 2.8285 $\pm$ 7.0927& 2.1544 $\pm$ 3.7220 \\
$1\,000$ &0.6490&  5911.86&  1631.04& 2.7014 $\pm$ 4.7655& 1.9614 $\pm$ 2.5031\\
$3\,000$ &0.6582&  4498.76&  0901.81& 2.4906 $\pm$ 2.4001& 1.7832 $\pm$ 1.0746\\
$5\,000$ &0.5998&  5239.97&  1542.03& 2.8828 $\pm$ 2.0064& 1.9242 $\pm$ 1.0885\\
$10\,000$ &0.6069&  5274.03&  1200.05& 2.7298 $\pm$ 1.4234& 1.8337 $\pm$ 0.6790\\
\bottomrule
\end{tabular}
\end{center}
\end{table}
\begin{table}[!htbp]
\caption{\label{SummurySCI3} Asymptotic variances and confidence
intervals ($u = u_{0.975}=1.96$) from subsamples of simulated under-varied trivariate data with $n=\mbox{10 000}.$}
\begin{center}
\begin{tabular}{rccccc}
\toprule
  $n\;$ & $\det\widehat{\boldsymbol{\rho}}$ & $\widehat{\sigma}_{gvi}^2$ & $\widehat{\sigma}_{mvi}^2$ & $\widehat{\mathrm{GVI}}_n\pm u\widehat{\sigma}_{gvi}/\!\sqrt{n}$ & $\widehat{\mathrm{MVI}}_n\pm u\widehat{\sigma}_{mvi}/\!\sqrt{n}$ \\
\toprule
50 & 0.9174&  200.3427& 180.0472& 0.7795 $\pm$ 3.9233& 0.8794 $\pm$ 3.7193 \\
100 &0.9634&  77.9392& 67.7172& 0.7354 $\pm$ 1.7303& 0.8242 $\pm$ 1.6129 \\
300 &0.9551&  76.8087& 69.3033& 0.6743 $\pm$ 0.9917& 0.7955 $\pm$ 0.9420 \\
500 &0.9446&  65.1680&  58.5276& 0.6174 $\pm$ 0.7076& 0.7309 $\pm$ 0.6706 \\
1 000 &0.9281&  49.6498&  44.1097& 0.5368 $\pm$ 0.4367& 0.6490 $\pm$ 0.4116\\
3 000 &0.9262&  34.0762&  29.2322& 0.4619 $\pm$ 0.2089& 0.5675 $\pm$ 0.1935\\
5 000 &0.9221&  32.6305&  28.0190& 0.4529 $\pm$ 0.1583& 0.5661 $\pm$ 0.1467\\
10 000 &0.9195 & 38.7897&  33.6378& 0.4980 $\pm$ 0.1221& 0.6161 $\pm$ 0.1137\\
\bottomrule
\end{tabular}
\end{center}
\end{table}
\begin{figure}[hbtp]
\begin{minipage}[c]{.415\linewidth}
\begin{center}
(a) $k=4$ over-: $\mathrm{GVI}=2.1905$
\includegraphics[scale=0.17]{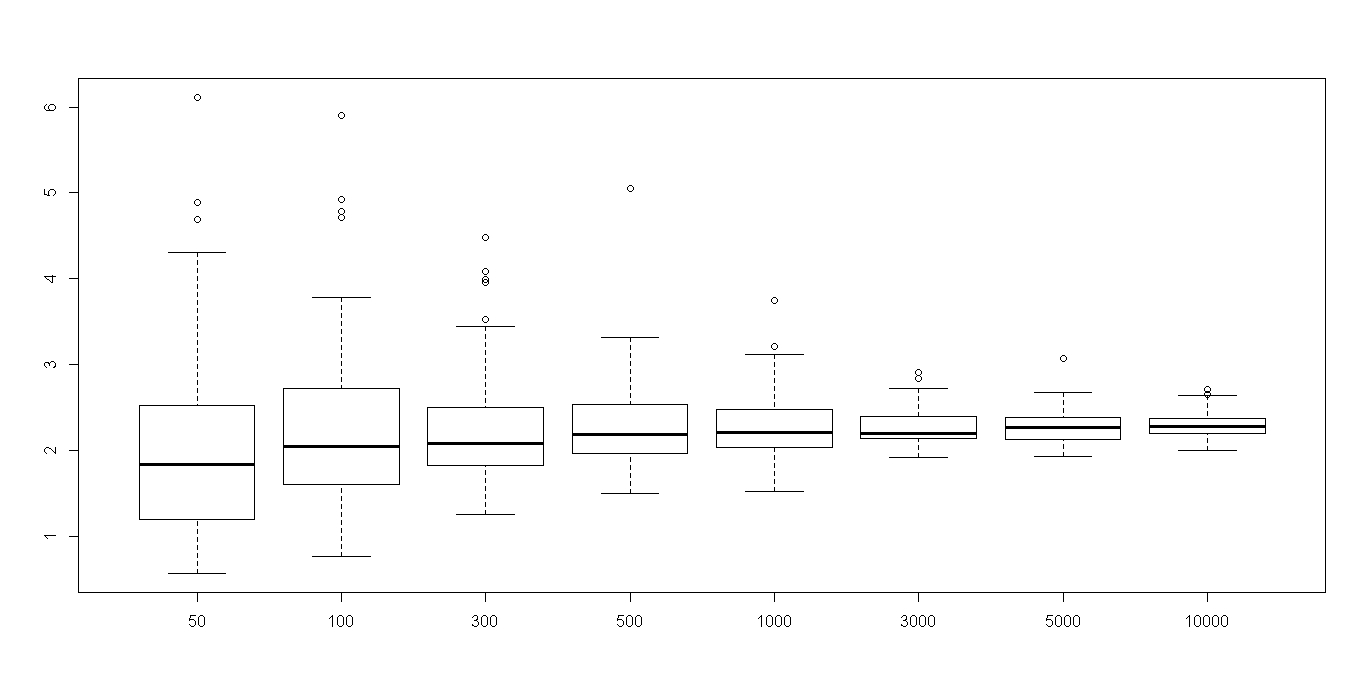}
\label{GDI}
\end{center}
\end{minipage}
\hfill
\begin{minipage}[c]{.415\linewidth}
\begin{center}
(b) $k=4$ over-: $\mathrm{MVI}=1.5827$
\includegraphics[scale=0.17]{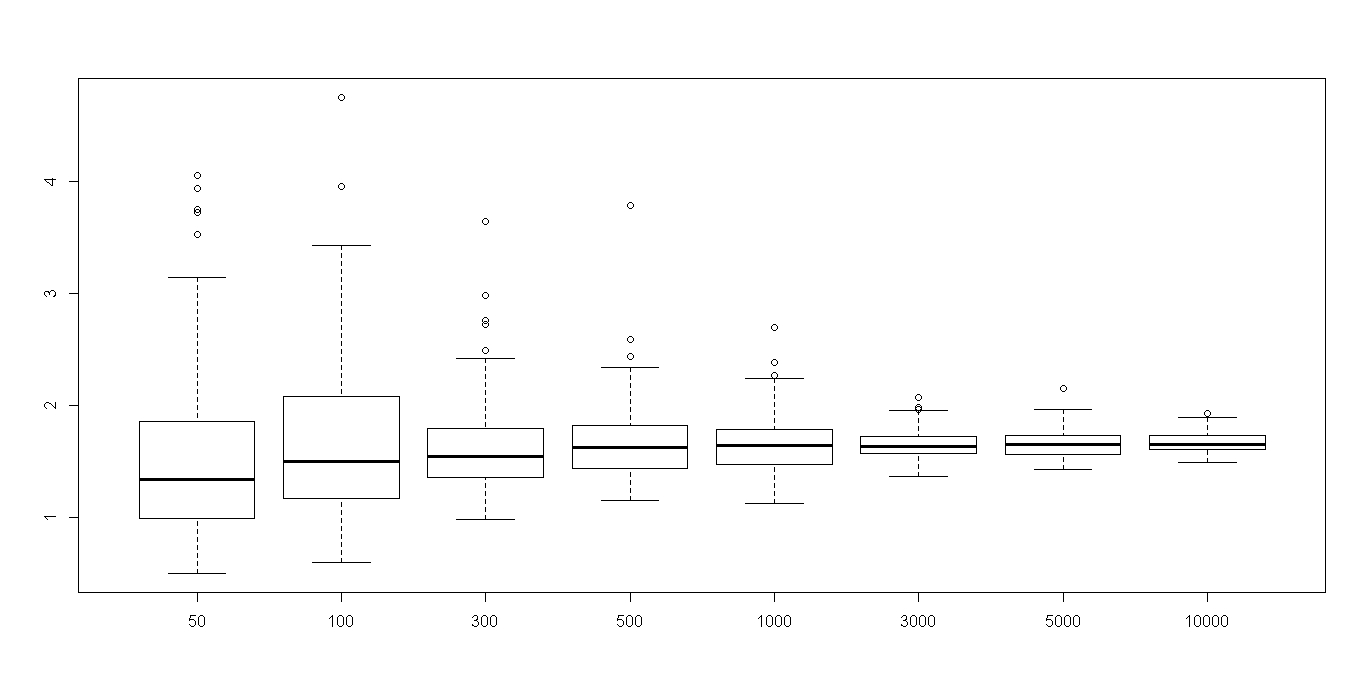}
\label{MDI}
\end{center}
\end{minipage}
\end{figure}
\begin{figure}[!hbtp]
\begin{minipage}[c]{.415\linewidth}
\begin{center}
(c) $k=6$ equi-: $\mathrm{GVI}=1.0259$
\includegraphics[scale=0.17]{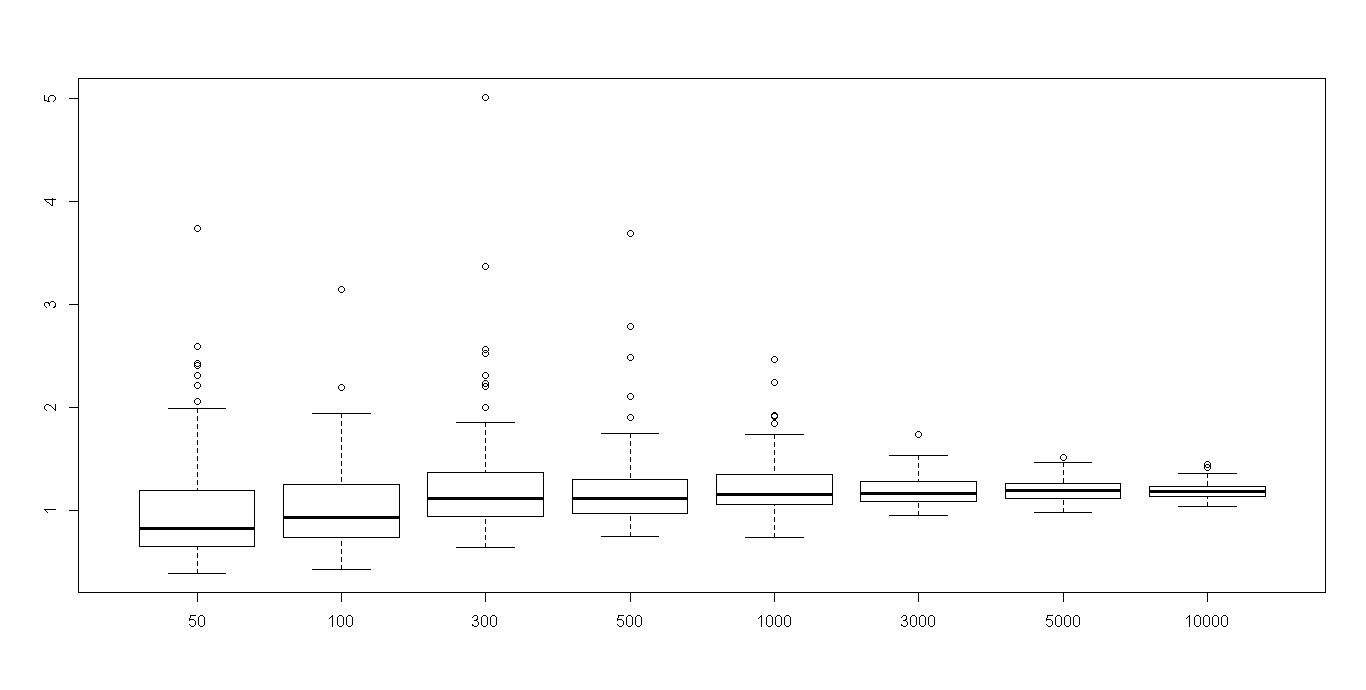}
\label{GDI}
\end{center}
\end{minipage}
\hfill
\begin{minipage}[c]{.415\linewidth}
\begin{center}
(d) $k=6$ equi-: $\mathrm{MVI}=0.9738$
\includegraphics[scale=0.17]{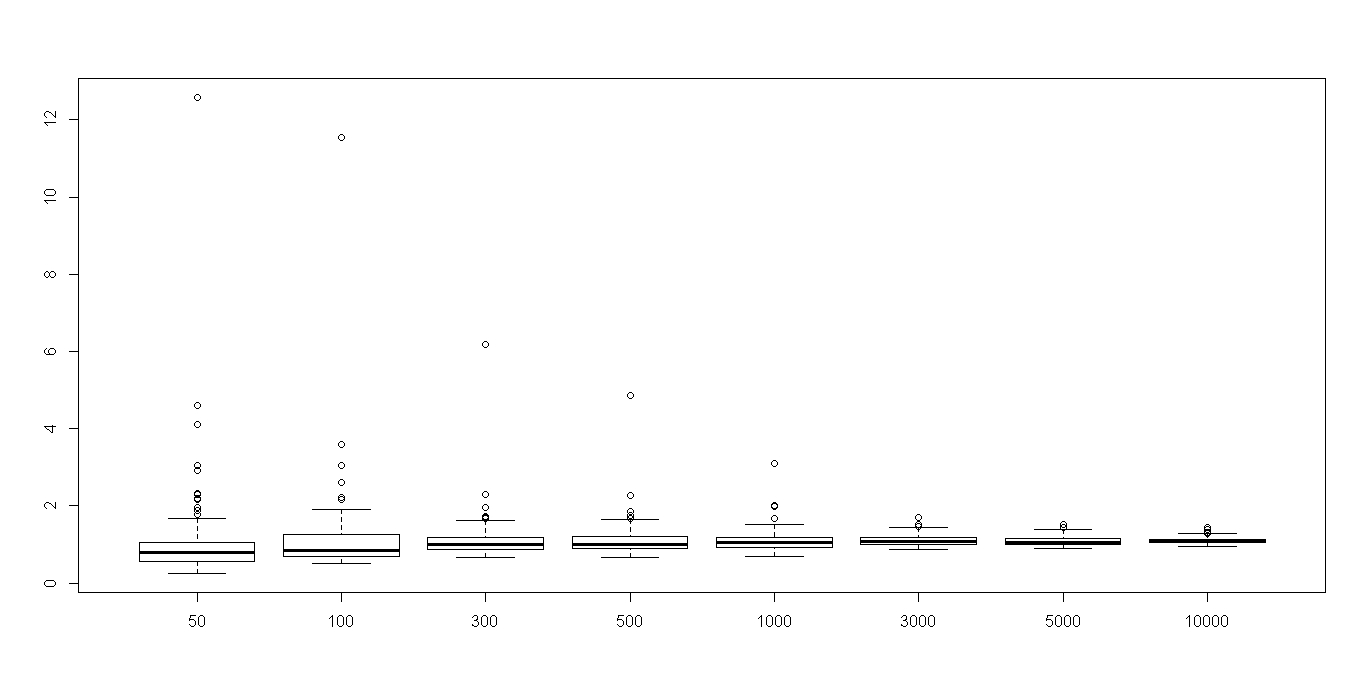}
\label{MDI}
\end{center}
\end{minipage}
\end{figure}
\begin{figure}[hbtp]
\begin{minipage}[c]{.415\linewidth}
\begin{center}
(e) $k=3$ under-: $\mathrm{GVI}=0.4425$
\includegraphics[scale=0.17]{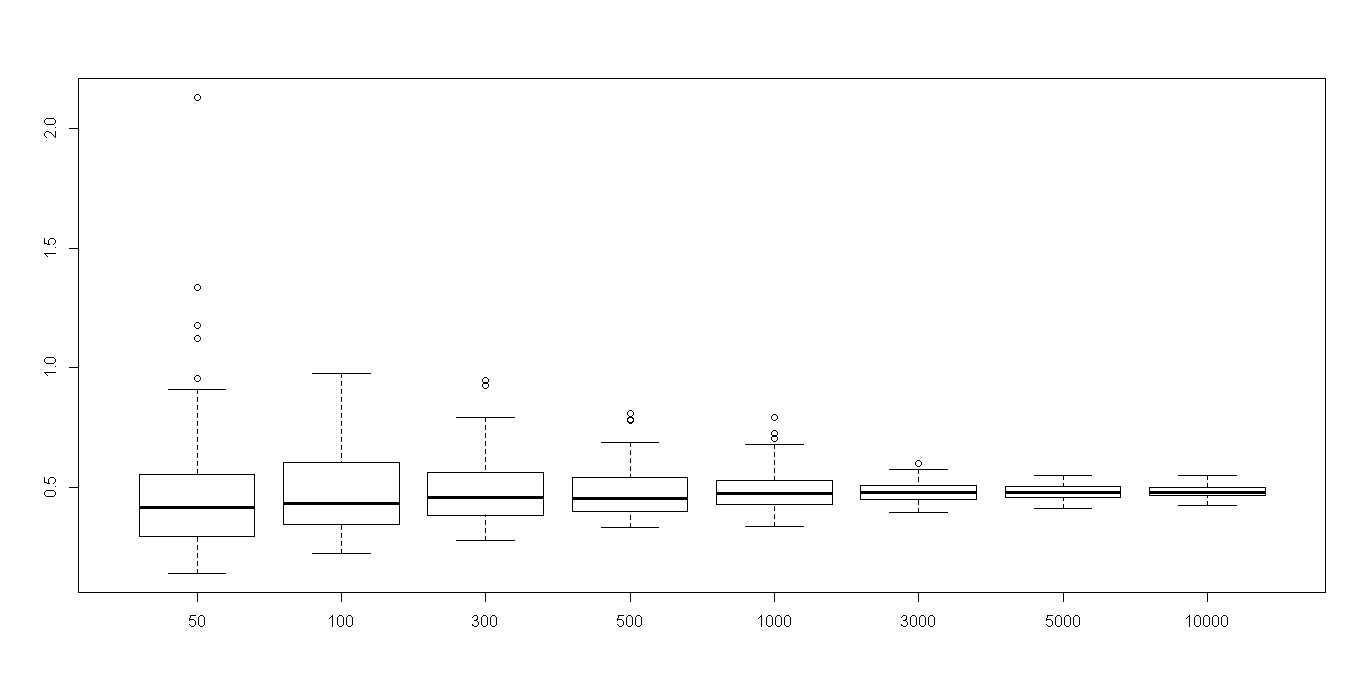}
\label{GDI}
\end{center}
\end{minipage}
\hfill
\begin{minipage}[c]{.41\linewidth}
\begin{center}
(f) $k=3$ under-: $\mathrm{MVI}=0.5635$
\includegraphics[scale=0.17]{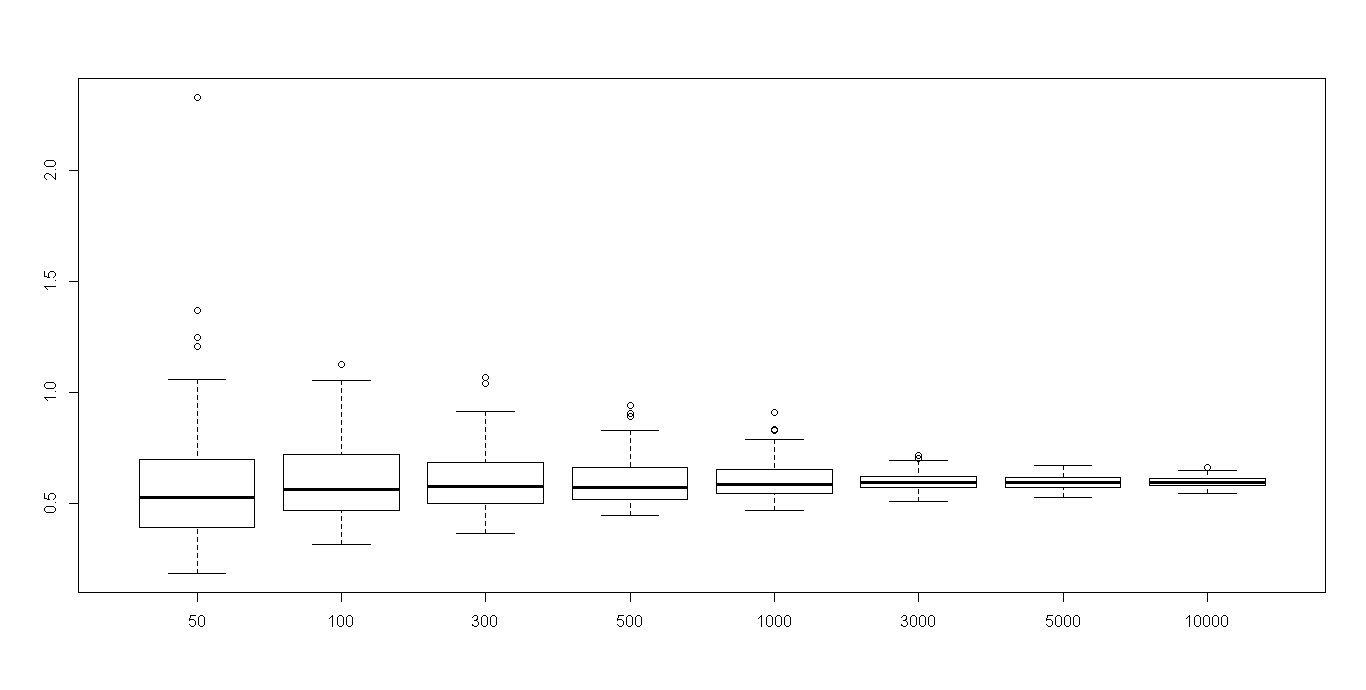}
\label{MDI}
\end{center}
\end{minipage}
\caption{Boxplots for the targets GVI and MVI with 100 replicates according to sample size and simulation parameters of Tables~\ref{SummurySCI}, \ref{SummurySCI2} and \ref{SummurySCI3}.}
\label{fig:boxpl}  
\end{figure}
Table~\ref{SummurySCI} depicts an evolution of the asymptotic variances and confidence intervals of $\mathrm{GVI}$ and $\mathrm{MVI}$ from subsamples of a simulated $6$-variate dataset with a maximum size $n=\mbox{10 000}$, having the same parameters as those for Table~\ref{SummurySsd}. We observe that both estimated standard errors $u\widehat{\sigma}/\!\sqrt{n}$ decrease when
sample size $n$ increases, and similarly in this context of $6$-variate quasi-equi-variation for GVI and also  for MVI. The stable behavior of the variances agrees with Proposition~\ref{Prop-DeltaGVI}.

Similar studies have also been performed simulating a $4$-variate over-varied distribution (Table~\ref{SummurySCI2}) and a trivariate under-varied distribution (Table~\ref{SummurySCI3}). The results shown in Table~\ref{SummurySCI2} have been obtained by simulating four marginal (over-, equi-, and under-varied) Weibull distributions, with the cross correlation matrix such that 
$$
\det\boldsymbol{\rho} := \det\left(
\begin{array}{cccc}
 1 & -0.03 & 0.27&0.35\\
 -0.03 & 1 & -0.05& 0.04   \\
   0.27 & -0.05 & 1& 0.48 \\
   0.35& 0.04& 0.48& 1\\
   \end{array}
\right)=0.6588.
$$
For the results shown in Table~\ref{SummurySCI3}, we have simulated one marginal Weibull distribution, one marginal exponential distribution and one marginal lognormal distribution, with the correlation matrix such that
$$
\det\boldsymbol{\rho} :=\det\left(
\begin{array}{cccc}
 1 & -0.03 & -0.27\\
 -0.03 & 1 & -0.05   \\
   -0.27 & -0.05 & 1\\
   \end{array}
\right)=0.9229.
$$
We also notice that all estimated standard errors $u\widehat{\sigma}/\!\sqrt{n}$ decrease when sample size $n$ increases, but more slowly for GVI than for MVI in Table~\ref{SummurySCI2} of the $4$-variate phenomenon of over-variation. Figure~\ref{fig:boxpl} clearly points out typical  behaviors of boxplots related to Tables~\ref{SummurySCI}, \ref{SummurySCI2} and \ref{SummurySCI3}. However, the estimated variances in Table~\ref{SummurySCI2} of the $4$-variate over-variation are much larger than those in Table~\ref{SummurySCI3} of the $3$-variate under-variation. Therefore, for small and moderate sample sizes one can use a bootstrapped approach or a robust version for reducing the estimated variances.
\begin{table}[!htbp]
\caption{\label{SummurySCI4}Asymptotic and bootstrap confidence intervals for GVI and MVI indexes from sample size $n$ simulated with $u = u_{0.975}=1.96$ using the dataset of Table~\ref{SummurySsd}.} 
\begin{center}
\begin{tabular}{rcccc}
\toprule
$n\;$ &  $\widehat{\mathrm{GVI}}_n\pm u\widehat{\sigma}_{gvi}/\!\sqrt{n}$ & Boots($\widehat{\mathrm{GVI}}\pm u\widehat{\sigma}/\!\sqrt{n}$)& $\widehat{\mathrm{MVI}}_n\pm u\widehat{\sigma}_{mvi}/\!\sqrt{n}$ & Boots($\widehat{\mathrm{MVI}}\pm u\widehat{\sigma}/\!\sqrt{n}$)\\
\toprule
30&1.0154 $\pm$ 38.1119&0.9603 $\pm$ 0.0869&0.9798 $\pm$ 37.9062&0.9656 $\pm$  0.0851\\
50&1.0110 $\pm$ 36.1508&0.9604 $\pm$  0.0498&1.0407 $\pm$ 36.0743&1.0149 $\pm$  0.0486\\
100& 1.0241 $\pm$ 26.1398&1.0119 $\pm$ 0.0359&0.9950 $\pm$ 26.0580&0.9620 $\pm$ 0.0352\\
300& 0.9715 $\pm$ 23.5589&1.0416 $\pm$ 0.0310&1.0703 $\pm$ 23.4215&1.0409 $\pm$ 0.0305\\
500& 1.1679 $\pm$ 15.6334&1.0229 $\pm$ 0.0172&1.1648 $\pm$ 15.5693&1.0150 $\pm$ 0.0169\\
1 000& 1.1994 $\pm$ 09.0242&1.0315 $\pm$ 0.0095 &1.1952 $\pm$ 08.9738&1.0278 $\pm$ 0.0093\\
\bottomrule
\end{tabular}
\end{center}
\end{table}

Table~\ref{SummurySCI4} presents behaviors of both  asymptotic and bootstrap confidence intervals for GVI and MVI in the situations of small and moderate sample sizes (e.g., Angelo and Brian, 2019). For these $6$-variate equi-varied datasets, we still observe that all estimated standard errors decrease when sample size $n$ increases, but more sharply and very weakly through the bootstrap method.

\section{Concluding remarks and extensions}\label{sect-Concl}

From the univariate case of variation index (Abid et al., 2019b) and the multivariate dispersion indexes for count models (Kokonendji and Puig, 2018), we have first introduced multivariate variation indexes GVI, MVI and RVI for continuous  distributions on non-negative orthant. All these proposed indexes are easy to handle from a theoretical and practical point of view. Unlike the intuitive marginal variation index MVI, the index GVI takes into account the correlations between variables. The ratio of two GVI provides the index RVI for changing the reference distribution of the measure of over-, equi- and under-variation in the multivariate framework. The interpretation and some properties of GVI and MVI are provided. Also, the asymptotic variances of GVI and MVI obtained from Proposition \ref{Prop-DeltaGVI} seem to provide large standard errors for small and moderate sample sizes; they can be improved, for instance, through a bootstrap method. An example of real data analysis is presented, helping to select an appropriate multivariate model.

Then, from $\mathrm{RVI}_{\boldsymbol{X}}(\boldsymbol{Y})$ given in (\ref{MRVI}) one exactly obtains its equivalent (i.e., relative dispersion index) $\mathrm{RDI}_{\boldsymbol{X}}(\boldsymbol{Y})$ for count models by changing the support $\mathbb{S}=\mathbb{N}^k$ of $\boldsymbol{X}$ and $\boldsymbol{Y}$ (Formula (9) of Kokonendji and Puig, 2018). Concerning a generalization of the basical GVI of (\ref{GVI-def}) which is also considered as a particular RVI with respect to to the uncorrelated exponential model, the recent univariate unification of  dipersion and variation indexes by Tour\'e et al. (2019) is used in the multivariate framework of natural exponential families as follows. 
Let $\boldsymbol{X}$ and $\boldsymbol{Y}$ be two random vectors on the same support $\mathbb{S}\subseteq\mathbb{R}^k$ and assume $\boldsymbol{m}:=\mathbb{E}\boldsymbol{X}=\mathbb{E}\boldsymbol{Y}$, $\boldsymbol{\Sigma}_{\boldsymbol{Y}}:=\mathrm{cov}\boldsymbol{Y}$ and $\mathbf{V}_{F_{\boldsymbol{X}}}(\boldsymbol{m}):=\mathrm{cov}(\boldsymbol{X})$ fixed, then the
\textit{relative variability index} of $\boldsymbol{Y}$ with respect to $\boldsymbol{X}$ can be defined as
\begin{equation*}\label{RWI}
\mathrm{RWI}_{\boldsymbol{X}}(\boldsymbol{Y}):=\mathrm{tr}[\boldsymbol{\Sigma}_{\boldsymbol{Y}}\mathbf{W}^+_{F_{\boldsymbol{X}}}(\boldsymbol{m})]\gtreqqless 1,
\end{equation*}
where $\mathbf{W}^+_{F_{\boldsymbol{X}}}(\boldsymbol{m})$ is the unique Moore-Penrose inverse of the associated matrix $\mathbf{W}_{F_{\boldsymbol{X}}}(\boldsymbol{m}):=[\mathbf{V}_{F_{\boldsymbol{X}}}(\boldsymbol{m})]^{1/2}[\mathbf{V}_{F_{\boldsymbol{X}}}(\boldsymbol{m})]^{\top/2}$ to $\mathbf{V}_{F_{\boldsymbol{X}}}(\boldsymbol{m})$; see Appendix B for GVI. 
Thus, we unify the construction of GDI and GVI by choosing $\mathbf{W}_{F_{\boldsymbol{X}}}(\boldsymbol{m})=\sqrt{\boldsymbol{m}}\sqrt{\boldsymbol{m}}^\top$ and $\mathbf{W}_{F_{\boldsymbol{X}}}(\boldsymbol{m})=\boldsymbol{m}\boldsymbol{m}^\top$, respectively. Note that one can consider $\mathbf{V}_{F_{\boldsymbol{X}}}(\boldsymbol{m})$ as a particular case of the MST variance function of Section \ref{ssect:MST}; but, it will be equivalent to the proposed GVI via RVI for supports $\mathbb{S}=[0,\infty)^k$ of distributions. Tests of hypothesis relying on the corresponding  estimators as test statistics with their asymptotic normality distributions should be deduced.

Finally, let us note the following problems which are in advanced discussion. Is it possible to characterize first the univariate over-/under-variation with respect to exponential distribution through the weighted exponential distribution as the count case by Kokonendji et al. (2008)? See also Kokonendji (2014) for some references. Therefore, how to investigate the multivariate connections to over-, equi- and under-variation through $\boldsymbol{m}\mapsto\mathrm{GVI}_{F}(\boldsymbol{m})$ or $\mathrm{MVI}_{F}(\boldsymbol{m})$? How, for instance, to discriminate some closed distributions from these indexes? See, e.g., Dey and Kundu (2009) for a univariate case. Statistical tests of these multivariate variation indexes can be produced in the direction of Aerts and Haesbroeck (2017); see also Feltz and Miller (1996).

\section*{Appendix A. On a broader multivariate exponential distribution}

According to Cuenin et al. (2016), taking $p=2$ in their multivariate Tweedie (1984) models of flexible dependence structure, another way to define a $k$-variate exponential distribution is given by
$\mathscr{E}_k(\boldsymbol{\Lambda})$. The $k\times k$ symmetric variation matrix
$\boldsymbol{\Lambda}= (\lambda_{ij} )_{i,j \in \{1, \ldots, k\}}$
is such that $\lambda_{ij}=\lambda_{ji}\geq 0$, the mean of the marginal exponential is $\lambda_{ii}>0$, and the nonnegative correlation terms satisfy
\begin{equation}\label{poisson_corr}
\rho_{ij} = \frac{\lambda_{ij}}{\sqrt{\lambda_{ii}\lambda_{jj}}}\in [0,\min\{R(i,j),R(j,i)\}),
\end{equation}
with
$R(i,j) = \sqrt{\lambda_{ii}/\lambda_{jj}} \, (1-\lambda_{ii}^{-1}\sum_{\ell\neq
i,j}\lambda_{i\ell} )\in (0,1)$. The construction of Cuenin et al. (2016) is perfectly defined having $k(k+1)/2$ parameters as in
$\mathscr{E}_k(\boldsymbol{\mu},\boldsymbol{\rho})$.  Furthermore, we attain the exact bounds of the  correlation terms in (\ref{poisson_corr}). The main fact is that Cuenin et al. (2016) pointed out the construction and simulation of the negative correlation structure from the positive one of (\ref{poisson_corr}) by using the inversion method.

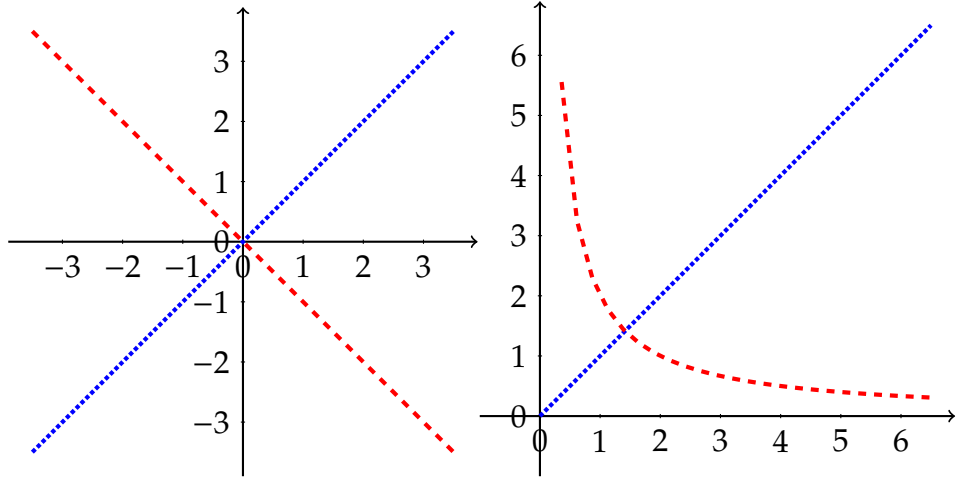
\begin{figure}[!htbp]
\begin{center}
\begin{tikzpicture}[scale=.8]
\draw[thick,->](-3.9,0)--(3.9,0);
\draw[thick,->](0,-3.9)--(0,3.9);
\foreach \x in {-3,-2,-1,0,1,2,3}
    \draw(\x,1pt)--(\x,-1pt)node[below]{$\x$};
\foreach \y in {-3,-2,-1,0,1,2,3}
    \draw(1pt,\y)--(-1pt,\y)node[left]{$\y$};
\draw[blue,ultra thick,densely dotted][domain=-3.5:3.5] plot(\x,\x);
\draw[red,ultra thick,dashed][domain=-3.5:3.5] plot(\x,-\x);
\end{tikzpicture}
\begin{tikzpicture}[scale=.8]
\draw[thick,->](-1,0)--(6.9,0);
\draw[thick,->](0,-1)--(0,6.9);
\foreach \x in {0,1,2,3,4,5,6}
    \draw(\x,1pt)--(\x,-1pt)node[below]{$\x$};
\foreach \y in {0,1,2,3,4,5,6}
    \draw(1pt,\y)--(-1pt,\y)node[left]{$\y$};
\draw[blue,ultra thick,densely dotted][domain=0:6.5] plot(\x,\x);
\draw[red,ultra thick,dashed][domain=0.36:6.5] plot(\x,2/\x);
\end{tikzpicture}
\caption{Support of bivariate distributions with maximum correlations (positive in blue and negative in red): model on $\mathbb{R}^2$ (left) and also finite support; model on $[0,\infty)^2$ (right), without finite support (Cuenin et al., 2016).}
\label{pict_corr}
\end{center}
\end{figure}

The negativity of a correlation component is important for the rare phenomenon of undervariation in a bivariate/multivariate positive continuous model. Figure~\ref{pict_corr} (right) plots
a limit shape of any bivariate positive continuous distribution with very strong negative correlation (in red), which is not the diagonal line of the upper bound ($+1$) of positive correlation (in blue); see, e.g., Cuenin et al. (2016) for bivariate count model. Contrarily, Figure~\ref{pict_corr} (left) represents the classic lower ($-1$) and upper ($+1$) bounds of correlations on $\mathbb{R}^2$ or finite support.

\section*{Appendix B. Construction of GVI}

In order to extend appropriately the univariate VI $=\sigma^2m^{-2}$ to the $k$-dimensional one for any positive continuous random vector $\boldsymbol{Y}$ on $(0,\infty)^k$ having positive (elementwise) mean vector $\boldsymbol{m}=(m_1,\ldots,m_k )^\top$ and covariance matrix $\boldsymbol{\Sigma}$, we consider the product of two matrices, namely $\boldsymbol{\Sigma}\boldsymbol{M}^{-1}$, where $\boldsymbol{M}=\boldsymbol{m}\boldsymbol{m}^\top$ is the $k\times k$ matrix outer product of $\boldsymbol{m}$ and which is well-defined. According to the singularity of $\boldsymbol{M}=\boldsymbol{m}\boldsymbol{m}^\top$, the unique Moore-Penrose inverse $\boldsymbol{M}^+$ of $\boldsymbol{M}$ is therefore 
$$\boldsymbol{M}^+=\frac{\boldsymbol{M}}{(\boldsymbol{m}^\top\boldsymbol{m})^2}.$$

Then, we have $\boldsymbol{\Sigma}\boldsymbol{M}^+=(\boldsymbol{M}^+\boldsymbol{\Sigma})^\top$. Since the rank of $\boldsymbol{M}$ is equal to 1, then $\boldsymbol{M}^+\boldsymbol{\Sigma}$ is also of rank 1 and has only one positive eigenvalue:
$$\lambda=\mathrm{tr}(\boldsymbol{M}^+\boldsymbol{\Sigma})=\mathrm{tr}(\boldsymbol{\Sigma}\boldsymbol{M}^+)=\frac{\boldsymbol{m}^\top\boldsymbol{\Sigma}\boldsymbol{m}}{(\boldsymbol{m}^\top\boldsymbol{m})^2}=:\mathrm{GVI},$$
where ``$\mathrm{tr}(\cdot)$'' stands for the trace operator.

This quantity $\lambda$ does not depend on the number $k$ of variables and it is numerically comparable to the univariate VI $=\sigma^2m^{-2}$. Also, it characterizes uniquely the $\boldsymbol{M}^+\boldsymbol{\Sigma}$ matrix, leading to the following definition of GVI. 
Note finally that if $\boldsymbol{\Sigma}=\boldsymbol{0}$ then we easily deduce $\lambda=0$, and conversely. 
We thus have the natural ordering of the half nonnegative real line for $\lambda\geq 0$.

\section*{Appendix C. Proofs of the asymptotic results}

\textbf{Proof of Proposition \ref{Prop-DeltaGVI}.} 
Part (i): Let
$\mathbf{Z}=(\ldots,Y_j,\ldots;\ldots,Y_jY_{\ell},\ldots)_{j\in \{1,\ldots,k\}; \ell \in \{j,\ldots,k\}}^\top$, 
$\mathbf{Z}_i=(\ldots,Y_{ij},\ldots;\ldots,Y_{ij}Y_{i\ell},\ldots)_{j\in \{1,\ldots,k\}; \ell \in \{j,\ldots,k\}}^\top$, 
for $i\in \{1,\ldots,n\}$, and the map $\Phi:(0,\infty)^k\times\mathbb{R}^{k(k+1)/2}\to (0,\infty)$ given through $\Phi(\mathbb{E}\mathbf{Z})=\mathrm{GVI}(\boldsymbol{Y})$ and $\Phi(n^{-1} \sum_{i=1}^n\mathbf{Z}_i)=\widehat{\mathrm{GVI}}_n(\boldsymbol{Y})$; i.e., for
$\boldsymbol{\theta}=(\ldots,m_j,\ldots;\ldots,\sigma_{j\ell},\ldots)_{j \in \{1,\ldots,k\}; \ell \in \{j,\ldots,k\} }^\top$,
 $\Phi(\boldsymbol{\theta})=({\boldsymbol{m}}^\top\boldsymbol{\Sigma}{\boldsymbol{m}})/(\boldsymbol{m}^\top\boldsymbol{m})^2$,
where $\boldsymbol{m}=(m_1,\ldots,m_k)^\top$ is the mean vector of $\boldsymbol{Y}$ and $\boldsymbol{\Sigma}= (\sigma_{j\ell} )_{j,\ell \in \{1,\ldots,k\}}$ is the covariance matrix of $\boldsymbol{Y}$ with $\sigma_{j\ell}=\sigma_{\ell j}$. Since $\Phi$ is differentiable at $\boldsymbol{\theta}$, the multivariate delta method (e.g., Serfling, 1980, Theorem A of Section 3.3) allows one to deduce that, as $n \to \infty$,
$$
\sqrt{n}\left\{\Phi\left(\frac{1}{n}\sum_{i=1}^n\mathbf{Z}_i\right)-\Phi(\mathbb{E}\mathbf{Z})\right\}\rightsquigarrow \mathcal{N}\left[ 0,\left( \frac{\partial\Phi(\boldsymbol{\theta})}{\partial\boldsymbol{\theta}}\right) ^\top\times \mathrm{cov}\boldsymbol{Z}\times \frac{\partial\Phi(\boldsymbol{\theta})}{\partial\boldsymbol{\theta}}\right].
$$

To check that $\mathrm{cov}\boldsymbol{\mathbf{Z}}=\boldsymbol{\Gamma}$ of the proposition  under the assumption on the fourth order moments of $Y_j$, one can rewrite $\mathbf{Z}$ in the following order:
$\mathbf{Z}=(\mathbf{Y},\widetilde{\mathbf{Y}})$ with $\boldsymbol{Y} = (Y_1,\ldots,Y_k)^{\top}$ and $\widetilde{\boldsymbol{Y}} = (\widetilde{Y}_1,\ldots,\widetilde{Y}_{\widetilde{k}})^{\top}$ such that $\widetilde{k}=k+(k-1)+\cdots+1=k(k+1)/2$ and 
\begin{multline*}
\widetilde{\boldsymbol{Y}} = (Y_1Y_1= \widetilde{Y}_1,Y_1Y_2= \widetilde{Y}_2,\ldots,Y_1Y_k= \widetilde{Y}_k, Y_2Y_2= \widetilde{Y}_{k+1},\ldots,Y_2Y_k= \widetilde{Y}_{2k-1},\\Y_3Y_3= \widetilde{Y}_{2k},\ldots,Y_kY_k= \widetilde{Y}_{k(k+1)/2})^{\top}. 
\end{multline*}
Then, the three main block matrices of $\boldsymbol{\Gamma}$ are successively found to be
\begin{align*}
\boldsymbol{\Sigma} &=\mathrm{cov}\boldsymbol{Y}= (\mathrm{cov}(Y_j,Y_{j'}) )_{j,j' \in \{1,\ldots,k\}}, \\
\boldsymbol{\Gamma}_4 & =
\mathrm{cov}\boldsymbol{\widetilde{Y}}= (\mathrm{cov}(\widetilde{Y}_j,\widetilde{Y}_{j'}))_{j,j' \in \{ 1,\ldots, k(k+1)/2\} }= (\mathrm{cov}(Y_{j_1}Y_{j_2}, Y_{j_1'}Y_{j_2'}))_{j_1,j_1' \in\{1,\ldots,k\};j_2 \in \{j_1,\ldots,k\} ;j_2' \in \{j_1',\ldots,k\}}, \\
\boldsymbol{\Gamma}_3^\top & = (\mathrm{cov}(Y_j,\widetilde{Y}_{j'}) )_{j \in \{1,\ldots,k\}; j' \in\{1,\ldots,   k(k+1)/2\}}= (\mathrm{cov}(Y_{j}, Y_{j_1'}Y_{j_2'}) )_{j,j_1' \in \{1,\ldots,k\} ;j_2'\in \{j_1',\ldots,k\}}.
\end{align*}

To see that $\partial\Phi(\boldsymbol{\theta})/\partial\boldsymbol{\theta}=\boldsymbol{\Delta}$, we first expand $\Phi$ as follows:
\begin{equation*}
\Phi(\boldsymbol{\theta})  =  \frac{\sum_{j=1}^k\sum_{\ell=1}^k\sigma_{j\ell}{m_jm_\ell}}{(\sum_{j=1}^k m_j^2)^2} =  \frac{ m_j^2\sigma_{jj}+{m_j}\sum\limits_{j'\neq j}\sigma_{jj'}{m_{j'}}+ \sum\limits_{\ell\neq j}{m_\ell} (\sigma_{\ell j}{m_j}+\sum\limits_{j'\neq j}\sigma_{\ell j'}{m_{j'}} )}{(m_j^2+\sum\limits_{j'\neq j}m_{j'}^2)^2}.
\end{equation*}
Then, direct calculations provide all components of $\boldsymbol{\Delta}$: for $j\in \{1,\ldots,k\}$, one has
\begin{eqnarray*}
\Delta_j=\frac{\partial}{\partial m_j}\, \Phi(\boldsymbol{\theta}) 
& = & \left\{2m_{j}\sigma_{jj}+\sum\limits_{j'\neq j}m_{j'}\sigma_{jj'}+\sum\limits_{\ell\neq j}m_{\ell}\sigma_{\ell j}-4m_j\left(\sum\limits_{j'= 1}^km_{j'}^2 \right) \,\Phi(\boldsymbol{\theta})\right\} (\boldsymbol{m}^\top\boldsymbol{m} )^{-2}\\
 & = & \left\{2\sum\limits_{j'= 1}^km_{j'}\sigma_{jj'}-4m_j\left(\sum\limits_{j'= 1}^km_{j'}^2 \right) \,\Phi(\boldsymbol{\theta})\right\} (\boldsymbol{m}^\top\boldsymbol{m} )^{-2},
\end{eqnarray*}
and $\Delta_{jj}=\partial\Phi(\boldsymbol{\theta})/\partial\sigma_{jj}=m_j^2/\left(\boldsymbol{m}^\top\boldsymbol{m}\right)^2$ while for $\ell \in \{j+1,\ldots, k\}$,
$\Delta_{j\ell}=\partial\Phi(\boldsymbol{\theta})/\partial\sigma_{j\ell}=2{m_jm_\ell}/\left(\boldsymbol{m}^\top\boldsymbol{m}\right)^2$. This ends the proof of Part (i).

\bigskip
\noindent
Part (ii): Introduce $\mathbf{W}=(Y_1,\ldots,Y_k;Y_1^2,\ldots, Y_k^2)^\top$, $\mathbf{W}_i=(Y_{i1},\ldots, Y_{1k};Y_{i1}^2,\ldots,Y_{ik}^2)^\top$ for $i \in \{ 1,\ldots,n\}$ and the map $\Psi:(0,\infty)^{2k}\to (0,\infty)$ defined by
$\Psi(\boldsymbol{\theta}) = (\sum_{j=1}^km_j^2\sigma_{jj})/(\sum_{j=1}^km_j^2)^2$
with $\boldsymbol{\theta}=(m_1,\ldots,m_k;\sigma_{11},\ldots,\sigma_{kk})^\top$. Then, one has $\Psi(\mathbb{E}\mathbf{W})=\mathrm{MVI}(\boldsymbol{Y})$ and $\Psi(n^{-1} \sum_{i=1}^n\mathbf{W}_i)=\widehat{\mathrm{MVI}}_n(\boldsymbol{Y})$. The function $\Psi$ is differentiable at the point $\boldsymbol{\theta}$ and, therefore, a straightforward application of the multivariate delta method leads to the conclusion that, as $n \to \infty$,
$$
\sqrt{n}\left\{\Psi\left(\frac{1}{n}\sum_{i=1}^n\mathbf{W}_i\right)-\Psi(\mathbb{E}\mathbf{W})\right\}\rightsquigarrow \mathcal{N}\left[0,\left( \frac{\partial\Psi(\boldsymbol{\theta})}{\partial\boldsymbol{\theta}}\right) ^\top\times \mathrm{cov}\boldsymbol{W}\times \frac{\partial\Psi(\boldsymbol{\theta})}{\partial\boldsymbol{\theta}}\right].
$$
Here, it is now trivial that $\mathrm{cov}\boldsymbol{W}=\boldsymbol{\Pi}$ of the theorem under the assumption of the finite moments on $Y_j$ and also that $\partial\Psi(\boldsymbol{\theta})/\partial\boldsymbol{\theta}=\boldsymbol{\Lambda}$ with 
$\Lambda_j=\partial\Psi(\boldsymbol{\theta})/\partial m_j =  \{2m_{j}\sigma_{jj}-4m_j\left( \sum_{j'=1}^km_{j'}^2\right)\,\Psi(\boldsymbol{\theta})\}\left(\boldsymbol{m}^\top\boldsymbol{m}\right)^{-2}$
and $\Lambda_{jj}=\partial\Psi(\boldsymbol{\theta})/\partial\sigma_{jj}=m_j^2\left(\boldsymbol{m}^\top\boldsymbol{m}\right)^{-2}$ for all $j\in \{1,\ldots,k\}$. This concludes the proof. \hfill $\blacksquare$

\bigskip
\textbf{Proof of Proposition \ref{Prop-sConsGVI}.}
According to the both continuous maps $\Phi:(0,\infty)^k\times\mathbb{R}^{k(k+1)/2}\to (0,\infty)$ defined through $\Phi(\mathbb{E}\mathbf{Z})=\mathrm{GVI}(\boldsymbol{Y})$ and $\Phi(n^{-1}\sum_{i=1}^n\mathbf{Z}_i)=\widehat{\mathrm{GVI}}_n(\boldsymbol{Y})$  and $\Psi:(0,\infty)^{2k}\to (0,\infty)$ such that $\Psi(\mathbb{E}\mathbf{W})=\mathrm{MVI}(\boldsymbol{Y})$ and $\Psi(n^{-1}\sum_{i=1}^n\mathbf{W}_i)=\widehat{\mathrm{MVI}}_n(\boldsymbol{Y})$ in the proof of Proposition \ref{Prop-DeltaGVI}, the desired result is easily deduced from
$n^{-1}\sum_{i=1}^n\mathbf{Z}_i \stackrel{a.s.}{\longrightarrow} \mathbb{E}\mathbf{Z}$ and $n^{-1}\sum_{i=1}^n\mathbf{W}_i \stackrel{a.s.}{\longrightarrow} \mathbb{E}\mathbf{W}$, respectively. \hfill $\blacksquare$

\end{document}